\newif\ifarXiv         
\newif\ifjournal        
\journalfalse            
\arXivtrue              
\documentclass[11pt]{article}
\usepackage[small,compact]{titlesec}
\usepackage{amsmath,amssymb,mathabx,setspace}
\usepackage{graphicx,epsfig,wrapfig,caption,epsfig,subcaption,sidecap}
\usepackage{url,color,verbatim,algorithmic}
\usepackage[ruled]{algorithm}  
\usepackage[sort,compress]{cite}
\usepackage[scaled]{helvet}
\usepackage[T1]{fontenc}
\usepackage[ansinew]{inputenc}
\usepackage[bookmarks=true, bookmarksnumbered=true, colorlinks=true,   pdfstartview=FitV,
linkcolor=blue, citecolor=blue, urlcolor=blue]{hyperref}
\usepackage[english]{babel}  
\usepackage{natbib} 
\usepackage{booktabs} 
\usepackage{MnSymbol}

\usepackage{multirow}
\let\oldbibliography\thebibliography
\renewcommand{\thebibliography}[1]{\oldbibliography{#1}
\setlength{\itemsep}{0pt}}

\def\wv{u}
\def\wu{ \widetilde{u}}  
\def\hatu{\widehat{u}}  
\def\hatf{\widehat{f}}    

\def\bA{\mathbf{A}}    
\def\bb{\mathbf{b}} 
\def\gap{\mathrm{Gap}}

\def\R{\mathbb{R}}       \def\C{\mathbb{C}}                               \def\E{\mathbb{E}}
                        
\newcommand{\Ebr}[1]{\mathbb{E}\left[{#1}\right]}

\newcommand{\red}[1]{\textcolor{black}{{#1}}}
\newcommand{\cyan}[1]{\textcolor{black}{{#1}}}

\newcommand{\fcommentout}[1]{}

\numberwithin{equation}{section}

\usepackage{soul}

\begin{document}

\begin{center}
\textbf{\Large Data-driven model reduction for stochastic Burgers equations } \\[0pt]
\vspace{4mm} Fei Lu\\
Department of Mathematics, Johns Hopkins University \\
 feilu@math.jhu.edu \\
\end{center}

We present a class of efficient parametric closure models for 1D stochastic Burgers equations. Casting it as statistical learning of the flow map, we derive the parametric form by representing the unresolved high wavenumber Fourier modes as functionals of the resolved variable's trajectory. The reduced models are nonlinear autoregression (NAR) time series models, with coefficients estimated from data by least squares. The NAR models can accurately reproduce the energy spectrum, the invariant densities, and the autocorrelations. 

Taking advantage of the simplicity of the NAR models, we investigate maximal space-time reduction. Reduction in space dimension is unlimited, and NAR models with two Fourier modes can perform well. The NAR model's stability limits time reduction, with a maximal time step smaller than that of the K-mode Galerkin system. We report a potential criterion for optimal space-time reduction: the NAR models achieve minimal relative error in the energy spectrum at the time step, where the K-mode Galerkin system's mean Courant--Friedrichs--Lewy (CFL) number agrees with that of the full model.
 
Key words: data-driven modeling, stochastic Burgers equation, closure model, CFL number. 
 
\tableofcontents


\section{Introduction}

 
Closure modeling aims for computationally efficiently reduced models for tasks requiring repeated simulations such as Bayesian uncertainty quantification \cite{stinisMoriZwanzigReduced2012,liIncorporationMemory2015} and data assimilation \cite{LTC17,LWM19}. Consisting of low-dimensional resolved variables, the closure model must take into account the non-negligible effects of unresolved variables \red{so as to capture both the short-time dynamics and large-time statistics}. As suggested by the Mori--Zwanzig formalism \cite{zwanzigNonequilibriumStatistical2001,CH13,LinLu20}, trajectory-wise approximation is no longer appropriate, and the approximation is in a statistical sense. That is, the reduced model aims to generate a process that approximates the target process in distribution, or at least, reproduce the key statistics and dynamics for the quantities of interest.  
For general nonlinear systems, such a reduced closure model is out of the reach of direct derivations from first principles. 

Data-driven approaches, which are based on statistical learning methods, provide useful and practical tools for model reduction. The past decades witness revolutionary developments of data-driven strategies, ranging from parametric models (see, e.g., \cite{kondrashov_Datadriven2015a,harlimParametricReduced2015,leiDatadrivenParameterization2016,xieDataDrivenFiltered2018,chekroun_Dataadaptive2017a,CL15,LLC17} and the references therein) to nonparametric and machine learning methods (see, e.g., \cite{pathakModelFreePrediction2018,ma2018model,harlim2020machine,parishParadigmDatadriven2016}). These developments demand a systematic understanding of model reduction from the perspectives of dynamical systems (see, e.g., \cite{duan2014effective,stinisRenormalizedMoriZwanzigreduced2015,LinLu20}), numerical approximation \cite{hudson2020coarse,choi2019space}, and statistical learning \cite{harlim2020machine,jiang2020_ModelingMissing}. 

With 1D stochastic Burgers equation as a prototype model, we aim to further the understanding of model reduction from an interpretable statistical inference perspective. More specifically, we consider a stochastic Burgers equation with a periodic solution on $[0, 2\pi]$: 
\begin{eqnarray} \label{SBE}
\begin{aligned}
&u_{t} = \nu u_{xx}-u u_{x}+ f(x,t) \,, 0 < x <2 \pi,t>0 \\
 &u(0,t)=u(2\pi,t) =0, \quad u_x(0,t)= u_x(2\pi,t), \\
\end{aligned}
\end{eqnarray}
from an initial condition $u(\cdot,0)$.  We consider a stochastic force $f(x,t)$ that is smooth in space, residing on $K_0$ low wavenumber Fourier modes, and white in time, given by 
\begin{equation} \label{eq:sforce}
 f(x,t)  = \sigma \sum_{m=1}^{K_0} \sin(mx) \dot{W}_m(t) + \cos(mx) \dot{W'}_m(t) ,
\end{equation}
where $\{W_m,W'_m\}$ are independent Brown motions. Here $\nu>0$ is the viscosity constant and $\sigma>0$ 
represents the strength of the stochastic force. 

Our goal is to find a discrete-time closure model for the first $K$ Fourier modes, so as to efficiently reproduce the energy spectrum and other statistics of these modes. 

We present a class of efficient parametric reduced closure models for 1D stochastic Burgers equations. The key idea is to approximate the discrete-in-time flow map statistically, in particular, to represent the unresolved high wavenumber Fourier modes as functionals of the resolved variable's trajectory. The reduced models are nonlinear autoregression (NAR) time series models, with coefficients estimated from data simply by least squares. We test the NAR models in four settings: reduction of deterministic responses ($K>K_0$) vs.~reduction involving unresolved stochastic force ($K<K_0$), and small vs.~large scales of stochastic force (with $\sigma=0.2$ and $\sigma=1$), where $K_0$ is the number of Fourier modes of the white-in-time stochastic force and $\sigma$ is the scale of the force. In all these settings, the NAR models can accurately reproduce the energy spectrum, invariant densities, and autocorrelation functions (ACF). We also discuss model selection, consistency of estimators, and memory length of the reduced models.

Taking advantage of our NAR models' simplicity, we further investigate a critical issue in model reduction of (stochastic) partial differential equations: maximal space-time reduction. The space dimension can be reduced arbitrarily in our parametric inference approach: NAR models with two Fourier modes perform well. The time reduction is another story.  The maximal time step is limited by the NAR model's stability and is smaller than those of the K-mode Galerkin system. Numerical tests indicate that the NAR models achieve the minimal relative error at the time step where the K-mode Galerkin system's mean CFL (Courant--Friedrichs--Lewy) number agrees with the full model's, suggesting a potential criterion for optimal space-time reduction. 

One can readily extend our parametric closure modeling strategy to general nonlinear dissipative systems beyond quadratic nonlinearities.  Along with \cite{LLC17}, we may view it as a parametric inference extension of the nonlinear Galerkin methods \cite{MT89,JKT90,Ros95,novo2001efficient}. However, it does not require the existence of an inertial manifold (and the stochastic Burgers equation does not satisfy the spectral gap condition that is sufficient for the existence of an inertial manifold \cite{Zel14}), and it applies to resolved variables of any dimension (e.g., lower than the dimension of the inertial manifold if it exists \cite{LLC17}). Notably, one may use NAR models that are linear in parameters and estimate them by least squares. Therefore, the algorithm is computationally efficient and is scalable for large systems. 

The limitation of the parametric modeling approach is its reliance on the derivation of a parametric form using the Picard iteration, which depends on the nonlinearity of the unresolved variables (see Section \ref{sec:derivParaM}). When the nonlinearity is complicated, a linear-in-parameter ansatz may be out of reach. One can overcome this limitation by nonparametric techniques \cite{jiang2020_ModelingMissing,zhang2019computing} and machine learning methods (see, e.g., \cite{harlim2020machine,pan2018data,ma2018model}). 
   
The stochastic Burgers equation is a prototype model for developing closure modeling techniques for turbulence (see e.g., \cite{e2000_InvariantMeasures,chorin_AveragingRenormalization2003,chorin2005_ViscositydependentInertial,bec2007_BurgersTurbulence,beck2009_UsingGlobal,wangTwolevelDiscretizations2011,dolaptchiev2013stochastic}). In particular, Dolaptchiev et al. \cite{dolaptchiev2013stochastic} propose a closure model for stochastic Burgers equation in a similar setting, based on local averages of finite-difference discretization, reproducing accurate energy spectrum similar to this study. We directly construct a simple yet effective NAR model for the Fourier modes, providing the ground of a statistical inference examination of model reduction. 
 
 \red{We note that the closure reduced models based on parametric inference are different from the widely studied proper orthogonal decomposition (POD)-based reduced order models (ROM) for parametric full models \cite{bennerSurveyProjectionBased2015,quarteroni2015reduced}. These POD-ROMs seek new effective bases to capture the effective dynamics by a linear system for the whole family of parametric full models. The inference-based closure models focus on nonlinear dynamics in a given basis and aim to capture both short-time dynamics and large-time statistics. In a probabilistic perspective, both approaches approximate the target stochastic process: the POD-ROMs are based on Karhunen-Lo\'eve expansion, while the inference-based closure models aim to learn the nonlinear flow-map. One may potentially combine the two and find nonlinear closure models for the nonlinear dynamics in the POD basis.}

The exposition of our study proceeds as follows. We first summarize the notations in Table \ref{tab:notations_all}.  Following a brief review of the basic properties of the stochastic Burgers equation and its numerical integration, we introduce in Section \ref{sec:SBE} the inference approach to closure modeling and compare it with the nonlinear Galerkin methods. Section \ref{sec3} presents the inference of NAR models: derivation of the parametric form, parameter estimation, and model selection.  Examining NAR models' performance in four settings in Section \ref{sec4}, we investigate the space-time reduction. Section \ref{sec5} concludes our main findings and possible future research.

\begin{table}[H]
\centering
\caption{Notations: the variables in the full and reduced models. }\label{tab:notations_all}
\begin{tabular}{c  l l } 
\toprule
      \textbf{Model }     & \textbf{Notation} & \textbf{Description}   \\
       \hline
 \multirow{6}{*}{Full model} &      $u(x,t)= \sum_{|k| \geq 1} \widehat u_k(t) e^{iq_kx}   $   & solution of \eqref{SBE} in its Fourier series   \\ \vspace{1mm}
       &   $f(x,t) = \sum_{|k| \geq 1}^{K_0} \widehat f_k(t) e^{iq_kx} $ &  stochastic force in \eqref{eq:sforce} in its Fourier series \\
      &  $v(x,t)=\sum_{|k| \leq K} \widehat u_k(t) e^{i q_k x} $  &  the resolved variable for closure modeling         \\
      &  $w(x,t)=\sum_{|k|> K} \widehat u_k(t) e^{i q_k x} $  &  the unresolved variable;  $u=v+w$ in \eqref{eq:w_intg} \\
       & $N$,$dt$ & number of modes and time step-size  \\
      \hline
\multirow{4}{*}{Reduced  models}    
      &  $K$ & number of modes in the NAR model \eqref{eq:NAR}  \\
      & $ (\wv_k^{n})_{|k|\leq K} $   & state variable in NAR, modeling $\widehat u_k(t_n)$ \\
       & $\delta = dt\times {\rm Gap}$  &  observation time interval \\  
      &  $R^\delta_k$, $\Phi^n$, $g^n$ & parametric terms in NAR, \eqref{eq:RM_discr} and \eqref{eq:NAR}  \\
\bottomrule
\end{tabular}
\end{table}

\section{Space-Time Reduction for Stochastic Burgers Equationations}\label{sec:SBE}
In this section, we first review basic properties of the stochastic Burgers equation and its numerical integration. Then, we introduce  inference-based model reduction and compare it with the nonlinear Galerkin methods.

\subsection{The Stochastic Burgers Equationation}
A Fourier transform of  Equation $\eqref{SBE}$ leads to 
\begin{eqnarray} 
\frac{d}{dt}\hatu_{k} 
&=&- \nu q_{k}^{2}\hatu_{k}- \frac{iq_{k}}{2}\sum_{l=-\infty}^{\infty}\hatu_{l}\hat{u}_{k-l}  + \widehat{f}_k(t) \label{FM} 
\end{eqnarray}
with $q_{k}=k, k \in \mathbb{Z}$, where $\hatu_{k}$ are Fourier modes: 
\begin{equation*}
\hatu_{k}(t)=\mathcal{F}[u]_{k}=\frac{1}{2\pi}%
\int_{0}^{2\pi}u(x,t)e^{-iq_{k}x}dx,\,\,\,u(x,t)=\mathcal{F}^{-1}[\hatu%
]=\sum_{k}\hatu_{k}(t)e^{iq_{k}x},
\end{equation*}%

The system has the following properties. First, it is Galilean invariant: if $u(x,t)$ is a solution, then $
u(x-ct,t)+c$, with $c$ an arbitrary constant speed, is a solution. To see this, let $v(x,t)=u(x-ct,t)+c$. Then, $v_t = -cu_x+ u_t$, $v_x=u_x$, and
\[v_t = cv_x + u_{xx}+ uu_x +f = cv_x + v_{xx} + (v-c)v_x +f = v_{xx} + vv_x +f. \] 
Without loss of generality, we set $\int_0^{2\pi} u(x,0)dx=0$. This implies that $\hatu%
_{0}(0)=0$. In this study, we only consider forces with mean zero, i.e. $\int_0^{2\pi} f(x,t) dx =0$, therefore from Eq.\eqref{FM}, we see that $\hatu_{0}(t)\equiv 0$, or equivalently, $\int_0^{2\pi} u(x,t)dx\equiv 0$. Second, the system has an invariant measure  \cite{sinai1991_TowResults,da2006introduction,e2000_InvariantMeasures}, due to a balance between the diffusion term, which dissipates energy, and the stochastic force, which injects energy.  In particular, the initial condition does not affect the large time statistical properties of the solution. Third, since $u$ is real, the Fourier modes satisfies $\hatu_{-k}=\hatu_{k}^{\ast }$, where $\hatu_{k}^{\ast }$ is the complex conjugate of $\hatu_{k}$. 

\subsection{Galerkin Spectral Method}\label{sec:num}
We consider the Galerkin spectral method for numerical solutions of the Burgers equation. The system is approximated as follows: the function 
$u(x,t)$ is represented at grid points $x_{i}=i \Delta x$ with $i=0,\dots ,2N-1$ and $\Delta x = \frac{2\pi}{2N}$. The Fourier
transform $\mathcal{F}$ is replaced by discrete Fourier transform 
\begin{equation*}
\hatu_{k}(t)=\mathcal{F}_{2N}[u]_{k}=\sum_{i=0}^{2N-1}u(x_{i},t)e^{-iq_{k}x_{i}},%
\,\,\,u(x_{i},t)=\mathcal{F}_{2N}^{-1}[\hatu]_{i}=\frac{1}{2N}%
\sum_{k=-N+1}^{N}\hatu_{k}e^{iq_{k}x_{i}}.
\end{equation*}%
For simplicity of notation, we abuse the notation  $u(x_{i},t)$ so that it denotes either the true solution or its high-resolution $2N$-mode approximation. 
Since $u$ is real, we have $\hatu_{-k}=\hatu_{k}^{\ast }$. Noticing
further that $\hatu_{0}=0$ due to Galilean invariance, and setting $\hatu%
_{N}=0$, we obtain a truncated system 
\begin{align}
\frac{d}{dt}\hatu_{k}=-\nu q_{k}^{2}\hatu_{k}-\frac{ik}{2} \sum_{\substack{|k-l|\leq N, |l|\leq N }}  \hatu_l \hatu_{k-l}+ \hat{f}_k, \text{ with } |k|=1,\dots ,N.  \label{DFM}
\end{align}

We solve Eq.\eqref{DFM} using the exponential time differencing fourth order Rouge--Kutta method (ETDRK4) (see \cite{CM02, KT05}) \cyan{with standard $3/2$ zero-padding for dealiasing (see e.g., \cite{GO77})}, with the force term $\widehat{f}_k$ treated as a constant in each time step. Such a mixture scheme \cyan{is of strong order 1, but it has an advantage of} preserving both the numerical stability of ETDRK4 and the simplicity of Euler--Maruyama. 

We will consider a relatively small viscosity $\nu=0.02$, so that random shocks are about to emerge in the solution. In general, a smaller viscosity constant demands a higher resolution in space-time to resolve the solution, particularly the emerging shocks as $\nu$ vanishes. To sufficiently resolve the solution, we set $N=128$ and $dt=0.001$. The solution is accurately resolved, with mean Courant--Friedrichs--Lewy (CFL) numbers being 0.139 and 0.045 for $\sigma=1$ and $\sigma=0.2$, respectively. Here the mean CFL number is computed as the average along a trajectory with $N_t=10^5$ steps 
\[ \text{Mean CFL number} =\frac{1}{N_t}\sum_{n=1}^{N_t} \sup_{x}|u(x,t_n)|\frac{\Delta t}{\Delta x},
\]
where $\Delta t$ and $\Delta x$ are the time step and space step, respectively. Furthermore, numerical tests show that the marginal densities converge as trajectory length increases.

\subsection{Nonlinear Galerkin and Inferential Model Reduction}
For simplicity of notation, we write the Burgers equation in an operator form as 
\begin{equation} \label{eq_diss} 
\partial_t u + Au = B(u) + f, \ \  u(0)= u_0
\end{equation}
with a linear operator $A:H^1_0(0,2\pi)\to L^2(0,2\pi)$ and a nonlinear operator $B: H^1_0(0,2\pi)\to L^2(0,2\pi)$
\[
A= -\nu \partial_{xx},\quad B(u) = -(u^2)_x/2,
\]
and with $f$ being the stochastic force. 

We first decompose the Fourier modes of $u$ into resolved and unresolved variables. 
 Recall that our goal of model reduction is to derive a closed system that can faithfully describe the dynamics of the coefficients $\{\widehat u_k(t)\}_{|k|=1}^K$, or equivalently, the low dimensional process $v( x,t)=\sum_{|k|=1}^K \widehat u_k(t) e^{iq_kx}$.

 Denote by $P$ the projection operator from $H_0^1(0,2\pi)$ to $\mathrm{span}\{e^{iq_kx}\}_{|k|=1}^K$, and let $Q:=I-P$ (and for simplicity of notation, we will also denote them as projections on the corresponding vector spaces of Fourier modes). With $u=Pu+ Qu=v+w$, 
 we can write the system \eqref{eq_diss} as 
 \begin{align}
\frac{dv}{dt}&=-PAv +PB(v) + Pf + [PB(v+w)-PB(v)],   \label{lowModes} \\
\frac{dw}{dt}&=-QAw+QB(v+w) + Qf.  \label{highModes}
\end{align}

To find a closed system for $v$, we quantify the truncation error  $PB(v+w)-PB(v)$ in \eqref{lowModes}, which represents the nonlinear interaction between the low and high wavenumber modes,  by either a function of $v$ or a functional of the trajectory of $v$.  In particular, in the nonlinear Galerkin method based on inertial manifold theory, see e.g.,  \cite{MT89, JKT90, Ros95,novo2001efficient}), one aims to represent the high modes $w$ as a function of the low modes $v$ (and hence obtaining an approximate inertial manifold). In the simplest implementation, one neglects the time derivative in Equationation \eqref{highModes} and solves $w=\psi(v)$ from
\begin{equation*}
w\approx -(QA)^{-1}[QB(v+w)+ Qf]
\end{equation*}
by fixed point iterations:
$\psi _{0}=0,\,\,\,\psi _{i+1}=-(QA)^{-1}[QB(u+\psi _{i}) + Qf]$. 
This leads to an approximation of $w$ as a function of $v$, which exists if $K$ is large enough and if the system satisfies a gap condition (so that an inertial manifold exists). However, among many dissipative systems with global attractor, only a few have been proven to satisfy the gap condition (see \cite{Zel14} for a recent review). More importantly, we can not always expect $K$ to be larger than the dimension of an inertial manifold, which is unknown in general. Therefore, such a nonlinear Galerkin approach works for neither a system without an inertial manifold nor for a $K$ smaller than the dimension of the inertial manifold. 

We take a different perspective on the reduction. Unlike the nonlinear Galerkin which aims for a trajectory-wise approximation, we aim for a probabilistic approximation of the distribution of the stochastic process $(v(\cdot, t),t\geq 0)$. The randomness of the process $v$ can come from random initial conditions and/or from the stochastic force. We emphasize that a key is to represent the dependence of the model error $PB(v+w)-PB(v)$ on the process $v$, not simply constructing a stochastic process with the same distribution as $PB(v+w)-PB(v)$, which may be independent of the process of $v$. 

In a data-driven approach, such a probabilistic approximation leads naturally to the statistical inference of the underlying process, aiming to represent the model error  $[PB(v+w)-PB(v)](t)$  as a functional of the past trajectory $(v(\cdot, s),s\leq t )$. This inferential reduction approach works flexibly for general settings: there is no need of an inertial manifold and the dimension $K$ can be arbitrary (e.g. less than the dimension of the inertial manifold, as shown in \cite{LLC17}). 

{\bf Space-time reduction.} To achieve a space-time reduction for practical computation, the reduced model should be a time series model with a time step $\delta>dt$ for time reduction, instead of a differential system. It approximates the flow map  (with $t_n=n\delta$)
\begin{equation}\label{eq:ukTrue}
\hatu_k(t_{n+1}) = F(\hatu_{\cdot}(t_n), \hatf_\cdot([t_n:t_{n+1}]))_k, \quad |k|\leq K,
\end{equation}
where $\hatu_{\cdot}(t_n) = (\hatu_{k}(t_n), |k|\geq 0)$ is the vector of all Fourier modes, and thus the above map is not a closed system for the low modes. Recall that for $|k|\leq K$, 
\begin{eqnarray} 
\underbrace{ \frac{d}{dt}\hatu_{k} 
=- \nu q_{k}^{2}\hatu_{k} -\frac{ik}{2} \sum_{\substack{|k-l|\leq K,\\ |l|\leq K }}  \hatu_l \hatu_{k-l}  }_{\text{$K$-mode truncation } } 
\underbrace{  - \frac{iq_{k}}{2}\sum_{\substack{|k-l|> K\\ \text{or }|l|> K }}\hatu_{l}\hatu_{k-l} }_{\text{truncation error}} + \widehat{f}_k(t)  \label{FM_decompose} 
\end{eqnarray}
Clearly, the K-mode truncated Galerkin system can provide an immediate approximation to $F$ in \eqref{eq:ukTrue}. Making use of it, we propose a time series model for $\{\hatu_k(t_n)\}_{|k|=1}^K$ in the form of 
\begin{equation}\label{eq:RM_discr}
 \wv_k^{n+1} =\wv_k^n + \delta [R^{\delta}_k(\wv^n) + f^n_k+ \Phi^n_k] + g^{n+1}_k , \quad |k|\leq K,
\end{equation}
where $R^{\delta}_\cdot(\wv^n)$ is from a one-step forward integrator with time step-size $\delta$ of the deterministic $K$-mode Galerkin, and $ f^n_k=\hatf_k(t_{n})$ is white noise in the $k$th Fourier mode of the stochastic force at time $t_n$. Here, the term $\Phi^n$ and the noise $g^{n+1}$ aim to represent the  truncation error, as well as the discretization error. Together with the other terms in \eqref{eq:RM_discr}, they provide a statistical approximation to the flow map $F$ in \eqref{eq:ukTrue}. 
In particular, the term $\Phi^n$ approximates $F$ based on information up to time $n$ (e.g., the conditional expectation),  and the noise $g^{n+1}$ aims to statistically represent the residual of the approximation.  Since the truncation error depends on the past history of the low wavenumber modes, and as suggested by the Mori--Zwanzig formalism \cite{CH13,LinLu20},  we make $\Phi^n$ depend on the trajectory $\wv^{1:n}$ of the state process, as well as the trajectories $f^{1:n}$ and $g^{1:n}$ : 
\begin{align}\label{eq:Phi_general}
 \Phi^n := \Phi(\wv^{1:n},f^{1:n}, g^{1:n}).
\end{align}
For simplicity, we assume the noise $\{g^n\}$ to be iid Gaussian, and the resulted time series model in \eqref{eq:RM_discr} is a nonlinear autoregression moving average model (NARMA) \cite{FY03,CL15,LLC16}. 

The right hand side of Eq.\eqref{eq:RM_discr}, together with  $\Phi^n$ defined in Eq.\eqref{eq:Phi_general},  aims for a statistical approximation of the discrete-time map \eqref{eq:ukTrue}. However, the general form in Eq.\eqref{eq:Phi_general} leads to a high dimensional function to be learned from data, which is intractable by regression methods using either global or local polynomial basis, due to the well-known curse of dimensionality. Fortunately, the physical model provides informative structures to reduce the dimension, and we can obtain effective approximations based on only a few basis functions with finite memory. In the next section, we derive from the physical model a parametric form for the reduced model, whose coefficients can be efficiently estimated from data.

To avoid confusions between notations, we summarize the correspondence of the variables between the full and reduced models in Table \ref{tab:FM-RMnotation}. 
\begin{table}[H]
\centering
\caption{Correspondence of the variables between the full and reduced models. }\label{tab:FM-RMnotation}
\begin{tabular}{l  c c} 
\toprule
          & \textbf{Full Model in} \eqref{DFM}   &  \textbf{NAR Model in}  \eqref{eq:RM_discr} \textbf{or} \eqref{eq:NAR} \\
          \hline
State variables   &   $ \hatu_k(t_n)$ or $\hatu (t_n)$ in \eqref{DFM} and \eqref{FM_decompose}   &  $ \wv^n_k$ or $u^n$  in \eqref{eq:RM_discr}  \\  
Resolved variable        & $v(x,t_n)$ or $v$, in \eqref{lowModes} and \eqref{eq:w_intg} &  vector $ (\wv^n_{-K},\ldots,\wv^n_{K})$ in \eqref{eq:NAR} \\
Unresolved variable        & $w(x,t)$ or $w$ in \eqref{highModes}  and \eqref{eq:w_intg}   & NA \\
 \hline 
\multirow{ 1}{*}{Stochastic force } & white noise $\hatf_k(t_ n) $ in  \eqref{FM_decompose}    & white noise $f^n_k$ in \eqref{eq:RM_discr}  \\
Noise introduced in NAR& NA   & $g^{n}$ in \eqref{eq:RM_discr} \\
\hline
Flow map &  $F$ in Eq.\eqref{eq:ukTrue} & Eq.\eqref{eq:RM_discr} \\
\bottomrule
\end{tabular}
\end{table}

\section{Inference of Reduced Models}
\label{sec3}
We present here the parametric inference of NAR models: derivation of parametric forms, estimation of the parameters, and model selection. 

\subsection{Derivation of Parametric Reduced Models}\label{sec:derivParaM}
We derive parametric reduced models by extracting basis functions from numerical integration of Eq.\eqref{lowModes}. The combination of these basis functions will give us $\Phi(\wv^{1:n},f^{1:n}, g^{1:n})$ in \eqref{eq:Phi_general}, which approximates the flow maps $\{ F(\hatu_{\cdot}(t_n), \hatf_\cdot([t_n:t_{n+1}]))_k, |k|\leq K\}$ in \eqref{eq:ukTrue} in a statistical sense. 

We first write a closed integro-differential system for the resolved process $(v(\cdot,t),t\geq 0)$. In view of Eq.\eqref{lowModes}, this can be simply done by integrating the equation of the high modes $w$ in  Eq.\eqref{highModes}: 
\begin{equation}\label{eq:w_intg}
\left \{
\begin{array}{ll}
\frac{dv}{dt}&=-PAv +PB(v) + Pf + [PB(v+w)-PB(v)],   \\
w(t) &= e^{-QA\tau} w(t-\tau) + \int_{t-\tau}^t e^{-QA(t-s)} [QB(v(s)+w(s))+ Qf(s)]ds,
\end{array}  
\right.
\end{equation}
where $\tau\in [0,t]$. Note that in addition to the trajectories $(v(\cdot, s), s\in [t-\tau,t] )$ and $(Qf(s), s\in [t-\tau, t])$, which we can assume to be known, the state $w(\cdot,t)$ also depends on the initial condition $w(\cdot, t-\tau)$. Therefore, this equation is not strictly closed. However, as $\tau$ increases, the effect of the initial condition decays exponentially, allowing for possible finite time approximate closure. Given  $w(\cdot, t-\tau)$ and  $(Qf(s), s\in [t-\tau, t])$,  the Picard iteration can provide us an approximation of  $w$ as a functional of the trajectory of $v$. That is, the sequence of functions $\{w^{(l)}\}$, defined by 
\begin{equation*}
w^{(l+1)}(t) = e^{-QA\tau} w^{(l)}(t-\tau) + \int_{t-\tau}^t e^{-QA(t-s)} [QB(v(s)+w^{(l)}(s))+ Qf(s)]ds,
\end{equation*}
with  $w^{(0)}(s)=0$ for $s\in [t-\tau, t]$, will converge to $w$ as $l\to \infty$. In particular, the first Picard iteration 
\begin{equation}\label{eq:w_intgAppr}
w^{(1)}(t) = \int_{t-\tau}^t e^{-QA(t-s)} [QB(v(s))+ Qf(s)]ds
\end{equation}
provides us a closed representation: from its numerical integrator, we can derive parametric terms for the reduced model. We emphasize that the goal is to derive parametric terms for statistical inference, but not to have a trajectory-wise approximation. Thus, high-order numerical integrators or high-order Picard iterations are helpful but may complicate the parametrization. For simplicity, we shall consider only the first Picard iteration and Riemann sum approximation of this integral.

We can now propose parametric numerical reduced models from the above integro-differential equation. In a simple form, we parametrize both the Riemann sum approximation of the first Picard iteration and a numerical scheme of the differential equation to obtain
\begin{align*} 
v (t_n)&\approx v(t_{n-1}) +a_1 \delta R^{\delta}(v(t_{n-1}) ) +a_2\delta Pf(t_{n-1}) + \delta[PB(v+w)-PB(v)](t_{n-1}),   \\
w(t_{n-1}) &\approx  \sum_{j=0}^p c_j e^{-QA j \delta} [ QB(v(t_{n-j})) +Qf(t_{n-j})] .
\end{align*}
Here $\delta = t_n-t_{n-1}$ denotes the time step-size, the nonlinear function $R^{\delta}(\cdot)$ comes from a numerical integration of the deterministic truncated Galerkin equation $\frac{dv}{dt} \approx -PAv +PB(v)$ at time $t_{n-1}$ and with time step-size $\delta$, and the coefficients $(a_1,a_2, c_j)$ are to be estimated by fitting to data in a statistical sense. To distinguish the approximate process in the reduced model from the original process, we denote it by $v^n$, and write the reduced model as
\begin{subequations}
 \begin{align} \label{eq:dis-sys}
v^n&=v^{n-1}+a_1 \delta R^{\delta}(v^{n-1} ) + a_2\delta Pf(t_{n-1}) + \delta[PB(v^{n-1}+w^{n-1})-PB(v^{n-1})] + g^n,   \\
w^{n-1} &= \sum_{j=1}^p c_j e^{-QA j \delta} [ QB(v^{n-j}) +Qf(t_{n-j})], \label{eq:QBvLag}
\end{align}
\end{subequations}
where $\{g^n\}$ is a process representing the residual, can be assumed to be stochastic force for simplicity, but can also be assumed to be a moving average part to better capture the time correlation as in \cite{CL15, LLC16}. The second Eq.\eqref{eq:QBvLag} does not have a residual term, as its goal is to provide a set of basis functions for the approximation of the forward map $v (t_n) = F(v(t_{n-1}), w(t_{n-1}), f)$ as in Eq.\eqref{eq:ukTrue}, not to model the high modes. 

Note that the time step-size $\delta$ can be relatively large, as long as the truncated Galerkin equation $\frac{dv}{dt} \approx -PAv +PB(v)$ of the slow variable $v$ can be reasonably resolved.  In general, such a step-size can be much larger than the time step-size needed to resolve the fast process $w$, because the effect of the unresolved fast process is ``averaged'' statistically when fitting the coefficients $(a_1,a_2, c_j)$ to data. Furthermore, the numerical error in the discretization is taken into account statistically.

Theoretically, the right-hand side of Eq.\eqref{eq:dis-sys} is an approximation of the conditional expectation $\Ebr{v(t_n)| v(t_{n-p:n-1}), Pf(t_{n-p:n-1})}$, which is the optimal $L^2$ estimator of the forward map conditional on the information up time $t_{n-1}$. Here, the $L^2$ is with respect to the joint measure of the vector $(v(t_{\cdot-p:\cdot-1}), Pf(t_{\cdot-p:\cdot-1}) )$, which is approximated by their joint empirical measure when fitting to data.

To avoid nonlinear optimization, the parametric form may be further simplified to be linearly dependent on the coefficients by dropping the terms that are nonlinear in the parameter, which is quadratic. In fact, recall that in the Burgers equation $B(u) = u u_x$ and  $PB(v+w)-PB(v) =v_x w  + vw_x + ww_x$. 
 By dropping the interaction between the high modes $ww_x$ and approximating 
 \[
 PB(v^{n-1}+w^{n-1})-PB(v^{n-1}) \approx v^{n-1}_x w^{n-1}  + v^{n-1}w^{n-1}_x
 \] 
in \eqref{eq:dis-sys}, we obtain a reduced model that depends linearly on the coefficients $\{a_j, c_j\}$.

\subsection{The Numerical Reduced Model in Fourier Modes} 

We now write the reduced model in terms of the Fourier modes as in Eq.\eqref{eq:RM_discr}. 

As discussed in the above section, the major task is to parametrize the truncation error $ PB(v+w)_k - PB(v)_k$. Recall that the operator $P$ projects $u$ to modes with wavenumber $1\leq |k| \leq K$ and that the bilinear function $PB(v)_k=\sum_{l} \hatu_l \hatu_{k-l}$ (hereafter, to simplify notation, we also denote $P$ and $Q$ on the corresponding vector spaces of Fourier modes). 
 \begin{equation}
 PB(v+w)_k - PB(v)_k = -\frac{ik}{2} \sum_{|l|>K \text{ or } |k-l|>K} \hatu_l \hatu_{k-l}.
 \end{equation}
Since the quadratic term $B(v)$ can only propagate energy from $(\hatu_{k}, 1\leq |k|\leq K)$ to modes with wave numbers less than $2K+1$,  we get only the high modes with wave numbers $K< |k|\leq 2K$ when we compute $w$ by a single iteration of $QB(v)$. (We use a single iteration for simplicity, but one can reach higher wave numbers by multiple iterations at the price of more complicated parametric forms.)
  Therefore, in a single iteration approximation, the truncated error will involve the first $2K$ Fourier modes: 
   \begin{equation*}
 PB(v+w)_k - PB(v)_k \approx -\frac{ik}{2} \sum_{\substack{ K<|k-l|\leq 2K \\ \text{ or } K<|l|\leq 2K }}  \hatu_l \hatu_{k-l}.
 \end{equation*}
 Dropping the interaction between the high-modes to avoid nonlinear optimization in parameter estimation, we have
  \begin{equation*}
 PB(v+w)_k - PB(v)_k \approx -\frac{ik}{2} \sum_{\substack{ |k-l|\leq K, K< |l| \leq 2K \\ \text{ or }  |l|\leq K, K< |k-l| \leq 2K} }   \hatu_l \hatu_{k-l}.
 \end{equation*}
We approximate the high modes $(\hatu_k, K< |k|\leq 2K)$ by  a functional of low modes as in \eqref{eq:QBvLag}, 
  \begin{equation*}
 \hatu_k(t_{n-1}) \approx  \sum_{j=1}^p c_{k,j} e^{-QAj\delta} [\wu_k(t_{n-j}) + \widehat{f}_k(t_{n-1})], \quad K< |k|\leq 2K
  \end{equation*}
 where $\wu_k$ is the high modes of the nonlinear function $B(v)$: 
 \[ \wu_k =QB(v)_k =-\frac{ik}{2}\sum_{|l|\leq K, |k-l|\leq K} \hatu_l \hatu_{k-l}, \text{ for } K< |k| \leq 2K.
 \]
  Here $QB(v)$ only represents the modes up to wavenumber $2K$, due to the fact that quadratic nonlinearity only involves interaction between double wave-numbers. One can reach higher wave numbers by iterations of the quadratic interaction.

The truncation error term can now be linearly parametrized as 
  \begin{align}\label{PBu-PBv}
 [PB(v+w) - PB(v)]_k(t_n)) &\approx -\frac{iq_k}{2} \sum_{j=0}^p  c_{k,j}  e^{-QAj\delta} \sum_{\substack{ |k-l|\leq K, K< |l| \leq 2K \\ \text{ or }  |l|\leq K, K< |k-l| \leq 2K} }\wu_l(t_{n})   \wu_{k-l}(t_{n-j}),
 \end{align}
 where we also denote $\wu_k= \hatu_k $ for $|k|\leq K$ for simplicity of notation. 
 \bigskip
 
 We have now reached a parametric numerical reduced model for the Fourier modes. Denote $\wv^n = (\wv^n_k, |k|\leq K)\in \C^K $ the low-modes in the reduced model that approximates the original low modes $(\hatu_k (t_n), |k|\leq K)$. The reduced model is 
 \begin{subequations} \label{eq:NAR}
 \begin{align} 
\wv^n_k & = \wv^{n-1}_k+  \delta [ R^{\delta}(\wv^{n-1}_\cdot) +   f^n_k + \Phi^n_k ]+ g^n_k, \quad 1\leq k\leq K, \label{eq:NAR_v} \\
\Phi^n_k & =  \sum_{j=1}^{p}\left[ c^v_{k,j}\wv^{n-j}_k + c^R_{k,j} R^{\delta}(\wv^{n-j}_\cdot) + c^f_{k,j} f^{n-j}_k+  c^w_{k,j} \sum_{\substack{ |k-l|\leq K, K< |l| \leq 2K \\ \text{ or }  |l|\leq K, K< |k-l| \leq 2K} }\wu^{n-1}_l \wu^{n-j}_{k-l} \right] \label{eq:NAR_linearQ} 
\end{align}
\end{subequations}
with the convention that $\wv^n_{-k}= (\wv^n_k)^*$ (with the sup-script $^*$ denoting
complex conjugate), and where the notion $\wu^{n-j}_l$  represents the high modes and is defined by  \begin{equation}  \label{eq:ks-ansatz-d}
  \widetilde{u}^{n-j}_{k }= 
    \left\{
    \begin{array}{ll}
      u^{n-j}_k~, & 1\leq k\leq K;\\[1ex]
      \frac{iq_k}{2}e^{-\nu q_k^2 j\delta} \sum_{|l|\leq K, |k-l|\leq K}\wv^{n-j}_{k-l} \wv^{n-j}_{l} , & K < k \leq 2K.
    \end{array}
    \right.
  \end{equation}
The reduced model is in the form of a nonlinear auto-regression moving average (NARMA) model: 

\begin{itemize}
\item The map $R^{\delta}(\cdot): \C^K \to \C^K$ is the 1-step forward of the deterministic $K$-mode Galerkin truncation equation $\frac{d v}{dt}=- PAv + PB(v) $ using a numerical integration scheme with a time step-size $\delta$, i.e. $v^{n+1}= v^{n} + \delta R^{\delta}(v^{n})$. 
We use the ETDRK4 scheme. 
\item The term $f^n_k$ denotes the increment of the $k$-th Fourier modes of the original stochastic force in the time interval $[t_{n-1},t_n]$, scaled by $1/\delta$, and it is separated from $R^{\delta}$ so that the reduced model can linearly quantify the response of the low-modes to the stochastic force. 

\item The function $\Phi^n_k :=\Phi^n_k(\wv^{n-p:n-1}, f^{n-p:n-1}) $ is a function $\C^{Kp+Kp}\to \C^K$   with parameters $\theta= (c^v, c^R,c^f,c^w)\in \R^{4Kp}$ to be estimated from data. In particular, the coefficients $c^v_{k,1}$ and $c^R_{k,1}$ act as a correction to the integration of the truncated equation. 

\item The new noise terms $\{g^n \in \C^K\}$ are assumed for simplicity to be a white noise independent of the original stochastic force $(f^n)$. That is, we assume that $\{g^n\}$ is a sequence of independent identically distributed (iid) Gaussian random vectors, with independent real and imaginary parts, distributed as $\mathcal{N}(0,\mathrm{Diag}(\sigma^g_k))$ with $\sigma^g_k$ to be estimated from data. Under such a white noise assumption, the parameters can be estimated simply by least squares (see next section).  In general, one can also assume other distributions for $g^n$, or other structures such as moving average $\{g^n:= \xi_n+ \sum_{j=1}^q c_j^g \xi_{n-j}\}$ with $\{\xi_n\}$ being a white noise sequence \cite{CL15,LLC16}. 
\end{itemize}


\subsection{Data Generation and Parameter Estimation} \label{sec:data_ParEst}
We estimate the parameters of the NAR model by maximizing the likelihood of the data. 

{\bf Data for the NAR model.} To infer a reduced model in form of Eq.\eqref{eq:NAR}, we generate relevant data from a numerical scheme that sufficiently resolve the system in space and time, as introduced in Section \ref{sec:num}. The relevant data are trajectories of the low-modes of the state and the stochastic force, i.e.  $\{\hatu_k(t_n), \widehat{f}_k(t_n)\}$ for $|k|\leq K$ and $n\geq 0$, which are taken as $\{\wv^n_k, f^n_k\}$ in the reduced model. Here, the time instants are $t_n= n\delta$, where $\delta$ can be much larger than the time step-size $dt$ needed to resolve the system. Furthermore, the data do not include the high modes. In short, the data are generated by a downsampling, in both space and time, of the high-resolution solutions of the system. 

The data can be either a long trajectory or many independent short trajectories. We denote the data consisting of $M$ independent trajectories by 
\begin{equation}\label{eq:data}
\text{Data: } \quad
\{\wv_k^{1:N_t,m}, f_k^{1:N_t,m}\}_{m,k=1}^{M,K} \text{ with } \wv_k^{1:N_t,m}= \hatu_k(t_{1:N_t})^{(m)}, f_k^{1:N_t,m} = \widehat{f}_k(t_{1:N_t})^{(m)} ,
\end{equation}
 where $m$ indexes the trajectories, $t_n = n\delta$ with $\delta$ being the time interval between two observations, and $N_t$ denotes the number of steps for each trajectory,

{\bf Parameter estimation.} The parameters in the discrete-time reduced model Eq.\eqref{eq:NAR} is estimated by maximum likelihood methods. Our discrete-time reduced model has a few attractive features: (i) the likelihood function can be computed exactly, avoiding possible approximation error that could lead to biases in estimators; (ii) the maximum likelihood estimator (MLE) may be computed by least squares under the assumption that the process $\{g^n\}$ is white noise, avoiding time-consuming nonlinear optimizations.

Under the assumption that $\{g^n\}$ is white noise, the parameters can be estimated simply by least squares, because the reduced model in Eq.\eqref{eq:NAR} depends linearly on the parameters. More precisely, the log-likelihood  of the data $\{\wv^{1:N_t,m}, f^{1:N_t,m}\}_{m=1}^M$ in \eqref{eq:data} can be written as
\begin{equation}\label{eq:lkhd}
l(\theta, \sigma^g) =- \sum_{|k|\leq K}\left[ \log\sigma^g_k + 
\sum_{n,m=1}^{T,M} \frac{ |\wv^{n,m}_k - \wv^{n-1,m}_k+   \delta  R^{\delta}(\wv^{n-1,m}_k) +\delta f^{n,m}_k + \delta\Phi^{n,m}_k (\theta)|^2}{2MT\sigma^g_k}  \right],
\end{equation}
where $|\cdot |$ denotes the absolute value of a complex number, $\theta= (c^v, c^R,c^f,c^w)\in \R^{4Kp}$ and $\sigma^g = (\sigma^g_1,\cdots, \sigma^g_K)\in \R^K$. To compute the maximum likelihood estimator (MLE) of the parameter $(\theta, \sigma^g)$, we note that  $\Phi^n_k (\theta) $ in \eqref{eq:NAR_linearQ} depends linearly on the parameter $\theta$. Therefore, the estimators of $\theta$ and $\sigma^g$ can be analytically computed by finding a zero of the gradient of the likelihood function. More precisely, 
denoting 
\[
\Phi_k^n(\theta) = \sum_{j=1}^{4p} \theta_j \Phi^n_{k,j}
\]
with $\Phi^n_{k,j}$ denoting the parameterized terms in \eqref{eq:NAR_linearQ},  we compute the MLE as
\begin{equation}\label{eq:MLE}
\begin{aligned}
\widehat \theta_k &= (\bA_k)^{-1}\bb_k, \quad 1\leq k \leq K, \\ 
\widehat \sigma^g_{k} & = \frac{1}{MT} \sum_{n,m=1}^{T,M} \|\wv_k^{n,m} - \wv^{n-1,m}_k+   \delta  R^{\delta}(\wv^{n-1,m}_k) + \delta f^{n,m}_k +  \delta \Phi^{n,m}_k (\widehat \theta)\|^2
\end{aligned}
\end{equation}
where the normal matrix $\bA_k$ and vector $\bb_k$ are defined by
\begin{align}
\bA_k(j',j) & =   \frac{\delta}{MT} \sum_{n,m=1}^{T,M} 	\langle \Phi^{n,m}_{k,j'}, \Phi^{n,m}_{k,j}\rangle, \quad 1\leq j', j \leq 4p,  \label{eq:normalMat}\\
\bb_k(j) & =    \frac{1}{MT} \sum_{n,m=1}^{T,M}	\langle \wv^{n,m}_k - \wv^{n-1,m}_k+   \delta  R^{\delta}(\wv^{n-1,m}_k) +\delta f^{n,m}_k , \Phi^{n,m}_{k,j}\rangle.  \notag
\end{align}
In practice, $A_k$ may be singular and it can be dealt with by pseudo inverse or regularization. 
We assume for simplicity that the stochastic force $g$ has independent components, so that the coefficients can be estimated by simple least square regression. One may further improve the NAR model by considering spatial correlation between the components of $g$ 
or by using moving average models \cite{CL15,LLC16,verheul2020stochastic} to account for the memory in the stochastic force. 


\subsection{Model Selection}\label{sec:modelSelection}
The parametric form in Eq.\eqref{eq:NAR_linearQ} leaves a family of reduced models with many freedoms underdetermined, such as the time lag $p$ and possible redundant terms. To avoid overfitting and redundancy, we proposed to select the reduced model by the following criterion. 
\begin{itemize}
\item Cross validation: the reduced model should be stable and can reproduce the distribution of the resolved process, particularly the main dynamical-statistical properties. We will consider  the energy spectrum, the marginal invariant densities, and temporal correlations:  
\begin{equation} \label{eq:stats} 
\begin{aligned} 
\text{Energy spectrum: } \E | \widehat u_k |^2 &= \lim_{N_tM\to\infty} \frac{1}{N_tM}\sum_{m,n=1}^{M,N_t}  | \widehat u_k(t_n)^{(m)} |^2; \\
\text{Invariant density of ${\rm Re} (\widehat u_k)$: } p_k(z)dz & =  \lim_{N_tM\to\infty} \frac{1}{N_tM}\sum_{m,n=1}^{M,N_t} \mathbf{1}_{ (z, z+ dz)}({\rm Re} (\widehat u_k(t_n)^{(m)}); \\ 
\text{Auto-correlation function: } \mathrm{ACF}_k(\tau) & = \E  [{\rm Re} \widehat u_k(t+\tau) {\rm Re} \widehat u_k(t) ] \\
& \approx \frac{1}{N_tM}\sum_{m,n=1}^{M,N_t} \mathrm{Re}(\widehat u_k(t_n+\tau)^{(m)}) \mathrm{Re}(\widehat u_k(t_n)^{(m)});
\end{aligned}
 \end{equation}
 for $k=1,\ldots,K$.

\item Consistency of the estimators. If the model is perfect and the data are either independent trajectories or a long trajectory from an ergodic measure, the estimators should converge as the data size increases (see e.g., \cite{FY03,Kut04}). While our parametric model may not be perfect, 
the estimators should also become less oscillatory as the data size increases, so that the algorithm is robust and can yield similar reduced models from different data sets. 
\item Simplicity and sparsity. When there are multiple reduced models performing similarly, we prefer the simplest model. We remove the redundant terms and enforce sparsity by LASSO (least absolute shrinkage and selection operator) regression \cite{tibshirani1996_RegressionShrinkage}. Particularly, a singular normal matrix \eqref{eq:normalMat} indicates the redundancy of the terms and the need to remove strongly correlated terms. 
\end{itemize}
These criteria are by no means exhaustive. Other methods include Bayesian information criterion (BIC, see, e.g., \cite{BD02}), and the error reduction ratio \cite{Bil13} may be applied, but in our experience, they provide limited help for the selection of reduced models \cite{LLC16,LLC17,LinLu20}.

In view of statistical learning of the high-dimensional nonlinear flow map  in \eqref{eq:ukTrue}, each linear-in-parameter reduced model provides an optimal approximation to the flow map in the function space spanned by the proposed terms. A possible future direction is to select adaptive-to-data hypothesis spaces in a nonparametric fashion \cite{jiang2020_ModelingMissing} and analyze the distance between the flow map and the hypothesis space spanned by these proposed terms \cite{Gyorfi06,LZTM19}.

\section{Numerical Study on Space-Time Reduction}
\label{sec4}
We examine the inference and performance of NAR models for the stochastic Burgers equation in \eqref{SBE} and \eqref{eq:sforce}. We will consider two settings of the full model: the stochastic force has a scale of either $\sigma =1$ or $\sigma =0.2$, representing that the stochastic force either dominates or subordinates to the dynamics, respectively. We will also consider two settings for reduction: the number of the Fourier modes of the reduced model is either $K>K_0$ or $K<K_0$, representing a reduction of the deterministic responses and a reduction involving stochastic force, respectively. 

\subsection{Settings}
As reviewed in Section \ref{sec:num}, we integrate the Eq.\eqref{DFM} of $2N$ Fourier modes by 
ETD-RK4 with a time-stepping $dt$ that the solution is resolved accurately. We call this discretized system the full model and its configuration is specified in Table \ref{tab:RMsettings}. We will consider two different scales for the stochastic force, with standard deviations $\sigma= 1$, leading to a dynamics dominated by the stochastic force, and $\sigma=0.2$, representing dynamics dominated by the deterministic drift. 
 \begin{table}[H]
  \caption{Settings of the full and reduced models}  \label{tab:RMsettings}
  \centering
\begin{tabular}{ c ll }
\toprule
\multirow{ 4}{*}{Full model} &  $\nu =0.02$,$L=1$ & viscosity, interval length of the equation   \\
                                            & $N=128$,$dt=0.001$ & number of modes, time step-size \\
                                           & $K_0=4$ & number of modes in the stochastic force\\
                                            & $\sigma = 1 \text{ or } 0.2$ & standard deviation of the stochastic force  \\
\hline
\multirow{3}{*}{Reduced  models} & $K=8 \text{ or } 2$     & number of modes in the reduced model\\ 
                                                    & $\delta = dt\times {\rm Gap}$  &  observation time interval\\  
                                                    & ${\rm Gap}\in \{5,10,20,30,40,50,80,160\}$ &  gap of time steps \\
\bottomrule
\end{tabular}
\end{table}

We generate data in \eqref{eq:data} from the full model as described in Section \ref{sec:data_ParEst}. We generate an ensemble of initial conditions by first integrating the system for $10^4$ time units from an initial condition $u_0(x)= \sin(x)+ 2 \cos(x)$ and draw $10^3$ samples uniformly from this long trajectory. Then, we generate either a long trajectory or an ensemble of trajectories starting from randomly picked initial conditions, and we save data with the time-stepping $\delta$.  Numerical tests show that the invariant densities and the correlation functions vary little when the data are generated from different initial conditions.

We then infer NAR models for the first $K$ Fourier modes with a time step $\delta$. We will consider two values for $K$ (recall that $K_0$ is the number of Fourier modes in the stochastic force)
\begin{itemize} 
\item $K=8>K_0=4$. In this case, $Qf=0$, i.e., the stochastic force does not act on the unresolved Fourier modes $w$ in \eqref{highModes}, so $w$ is a deterministic functional of the history of the resolved Fourier modes. In view of \eqref{eq:QBvLag}, the reduced model mainly quantifies this deterministic map. We call this case ``reduction of the deterministic response'' and present the results in Section \ref{sec:K=8}. 
\item $K=2<K_0$. In this case, $Qf \neq 0$, and $w$  in \eqref{highModes} depend on the unobserved Fourier modes of the stochastic force. Thus, the reduced model has to quantify the effects of the unresolved Fourier modes of both the solution and the stochastic force. We call this case ``reduction involving unresolved stochastic force'' and present the results in Section \ref{sec:K=2}. 
\end{itemize}
In either case, we explore the maximal time step that NAR models can reach by testing time steps $\delta = dt\times \{5,10,20,30,40,50,80,160\}$. 

We summarize the configurations and notations in Table \ref{tab:RMsettings}.

\subsection{Model Selection and Memory Length}\label{sec:model_selection}
We demonstrate model selection and the effect of memory length for reduced models with time step  $\delta=5dt$. We aim to select a universal parametric form of the NAR model for different setting of $(K,\sigma)$, where $K\in \{ 8,2\} $ is the number of Fourier modes in the NAR model and $\sigma\in \{1, 0.2\}$ is the standard deviation of the full model's stochastic force. Such a parametric form will be used later for the exploration of maximal time reduction by NAR models in the next sections. 

We select the model according to Section \ref{sec:modelSelection}: for each pair $(K,\sigma)$, we test a pool of NAR models and select the simplest model that best reproduces the statistics and has consistent estimators. The statistics are computed along a long trajectory of $T=2000$ time units. We say that an NAR is numerically unstable if it blows up (e.g. $|u^n|$ exceeding $10^5$) before reaching $T=2000$ time units.

We estimate the coefficients in \eqref{eq:NAR_linearQ} for a few time lag $p$s. Numerical tests show that the normal matrix in regression is almost singular, either when the stochastic force $f^{n-j}_k$ presents or when the lag for $u^{n-j}_k$ or $R^\delta(u^{n-j})$ is bigger than two. Thus, for simplicity, we remove them by setting: 
\begin{equation}\label{coefsTBE}
\begin{aligned}
 c^f_{k,j} = 0 \text{ for all } 1\leq j\leq p, \quad \text{and } c^{v}_{k,j} = c^{R}_{k,j}=0 \text{ for all } 1<j \leq p, 
\end{aligned}   
\end{equation}
and estimate only $c^{v}_{k,1}, \, c^{R}_{k,1}, \,  c^{w}_{k,j}$ for $1\leq j\leq p$. 

That is, in \eqref{eq:NAR_linearQ}, the terms $u^{n-j}_k$ and $R^\delta(u^{n-j})$ have a time lag 1, the stochastic force term $f^{n-j}_k$ is removed, and only the high-order (the fourth) term has a time lag $p$.  The memory length is $p\delta$.  

\begin{figure}[H]
    \centering
    \hspace{1mm}
  \ifarXiv
     \includegraphics[width=0.90\textwidth]{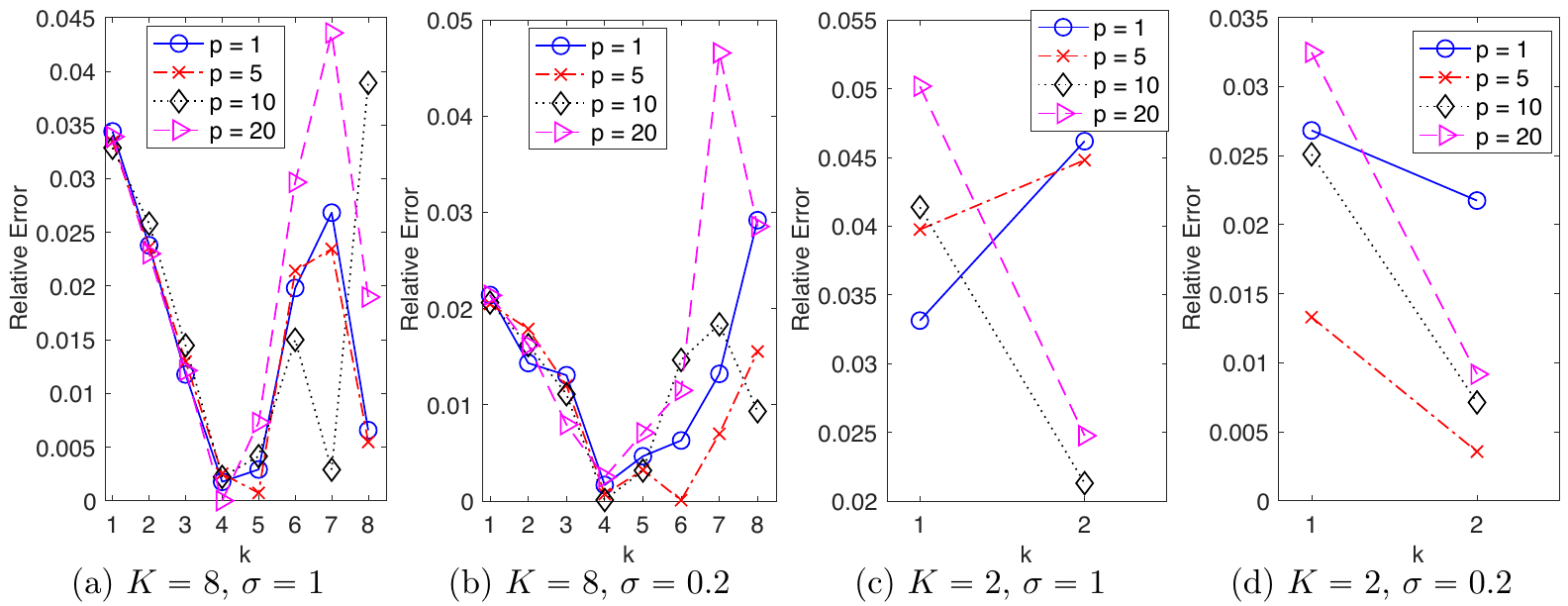} \\
\fi
 \ifjournal   
     \begin{subfigure}[b]{0.225\textwidth}
           \includegraphics[width=1\textwidth]{./figures/VspectrumL_gap5_K8f1.pdf}\vspace{-2mm}
         \caption{$K=8$, $\sigma=1$}
    \end{subfigure}
         \begin{subfigure}[b]{0.225\textwidth}
           \includegraphics[width=1\textwidth]{./figures/VspectrumL_gap5_K8fp2.pdf}\vspace{-2mm}
         \caption{$K=8$, $\sigma=0.2$}
    \end{subfigure}
         \begin{subfigure}[b]{0.225\textwidth}
           \includegraphics[width=1\textwidth]{./figures/VspectrumL_gap5_K2f1.pdf}\vspace{-2mm}
         \caption{$K=2$,  $\sigma=1$}
    \end{subfigure}
         \begin{subfigure}[b]{0.225\textwidth}
           \includegraphics[width=1\textwidth]{./figures/VspectrumL_gap5_K2fp2.pdf}\vspace{-2mm}
         \caption{$K=2$,  $\sigma=0.2$}
    \end{subfigure}
\fi
       \caption{Relative error in energy spectrum reproduced by the NAR models with different memory lengths $p$, in four settings of $(K,\sigma)$. As the time lag $p$ increases, the relative error tends to first decrease and then increase, particularly in (\textbf{b}) and (\textbf{d}) with $\sigma=0.2$. } \label{fig:spectrumLag}\vspace{-3mm}
\end{figure}

\vspace{2mm}
{\bf Memory length.} To select a memory length, we test NAR models with time lags $p\in \{ 1,5,10,20\}$ and consider their reproduction of the energy spectrum  in \eqref{eq:stats}. Figure \ref{fig:spectrumLag} shows the relative error in energy spectrum of these NAR models. It shows that as $p$ increases: (1) when the scale of the stochastic force is large ($\sigma=1$), the error oscillates without a clear pattern; (2) when $\sigma=0.2$, the error first decreases and then increases. Thus, a longer memory does not necessarily lead to a better reduced model when the  stochastic force dominates the dynamics; but when deterministic flow dominates the dynamics, a proper memory can be helpful. 

\vspace{-3mm}

\vspace{-3mm}
\begin{figure}[H]
    \centering
    \hspace{1mm}
 \ifarXiv
     \includegraphics[width=0.90\textwidth]{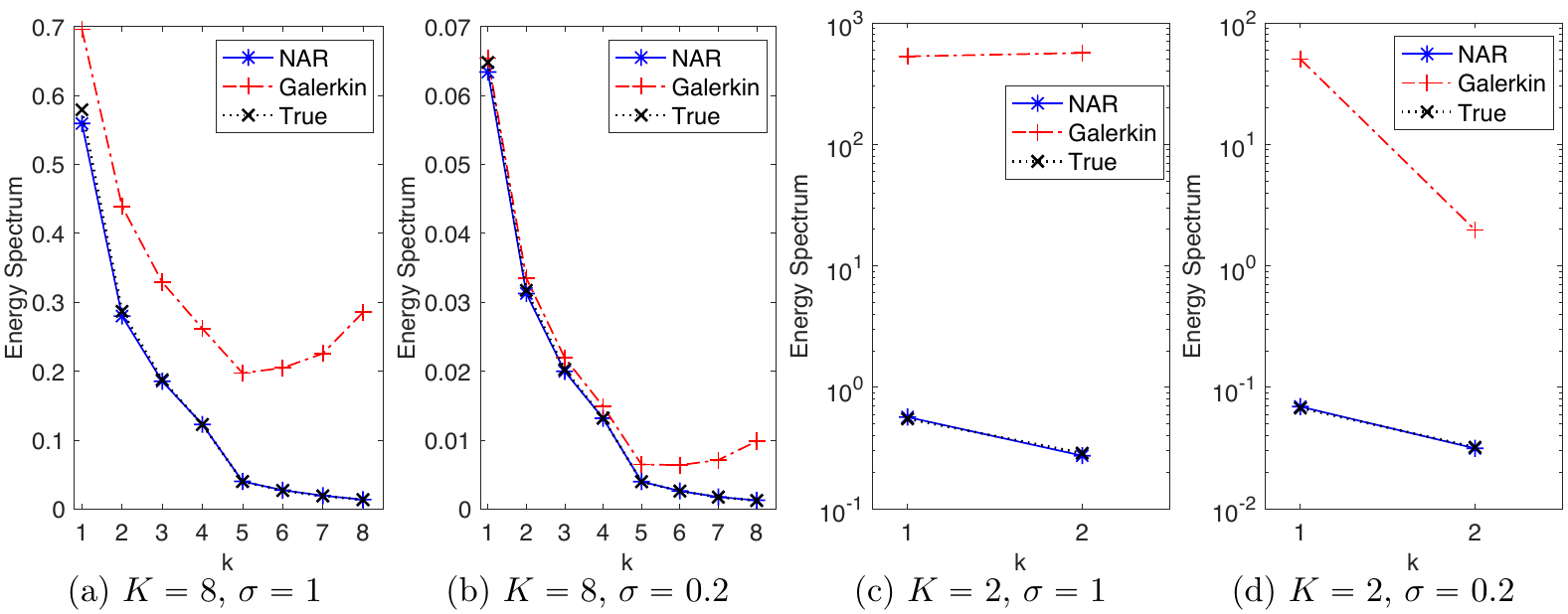} \\
\fi
 \ifjournal
     \begin{subfigure}[b]{0.225\textwidth}
           \includegraphics[width=1\textwidth]{./figures/VspectK8f1Gap5.pdf}\vspace{-2mm}
         \caption{$K=8$, $\sigma=1$}
    \end{subfigure}
         \begin{subfigure}[b]{0.225\textwidth}
           \includegraphics[width=1\textwidth]{./figures/VspectK8fp2Gap5.pdf}\vspace{-2mm}
         \caption{$K=8$, $\sigma=0.2$}
    \end{subfigure}
         \begin{subfigure}[b]{0.225\textwidth}
           \includegraphics[width=1\textwidth]{./figures/VspectK2f1Gap5.pdf}\vspace{-2mm}
         \caption{$K=2$,  $\sigma=1$}
    \end{subfigure}
         \begin{subfigure}[b]{0.225\textwidth}
           \includegraphics[width=1\textwidth]{./figures/VspectK2fp2Gap5.pdf}\vspace{-2mm}
         \caption{$K=2$,  $\sigma=0.2$}
    \end{subfigure}
 \fi
        \caption{Energy spectrum of NAR models with $p=1$ and the $K$-mode Galerkin systems in four settings of $(K,\sigma)$. The time step is $\delta = 5dt$ for the NAR models and is $dt$ for the Galerkin models. The NAR models accurately reproduce the true energy spectrum in all settings. } \label{fig:spectrum}
\end{figure}

In all four settings, the simplest NAR models with $p=1$ can consistently reproduce the energy spectrum with relative errors within 5\%. Remarkably,  the accuracy remains  when the true energy spectrum is at the scale of $10^{-2}$ for the modes with $k=7,8$ in Figure \ref{fig:spectrum}a,b and $k=2$ in Figure \ref{fig:spectrum}d. Figure \ref{fig:spectrum} also shows that the truncated $K$-mode Galerkin systems cannot reproduce the true energy spectrum in any of these settings, with upward tails, due to the lack of fast energy dissipation from the high modes. Thus, the NAR model has introduced additional energy dissipation through $\Phi^n$.

\begin{figure}[htb]    
    \centering     \hspace{1mm}
\ifarXiv
     \includegraphics[width=0.96\textwidth]{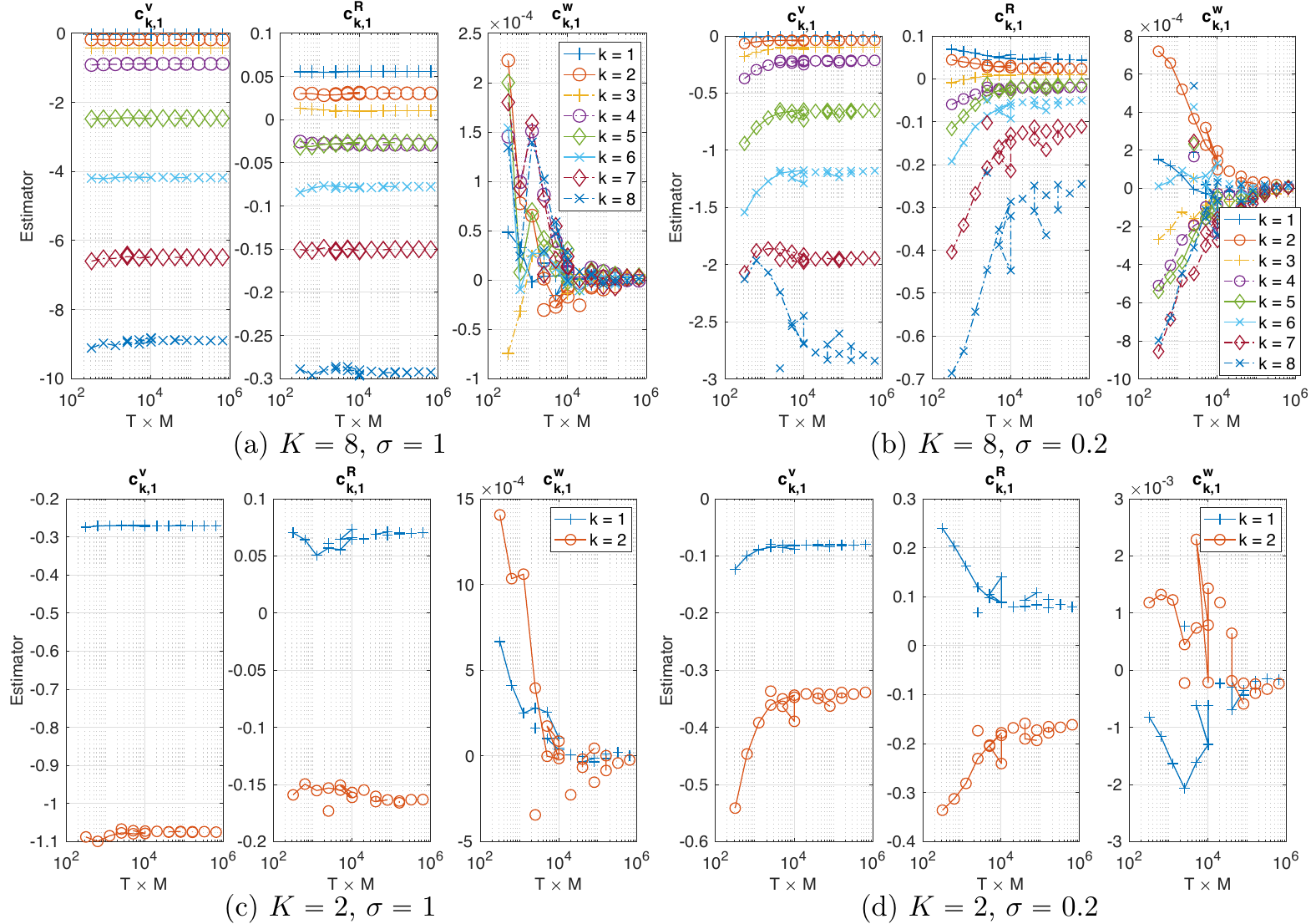} 
\fi
 \ifjournal
     \begin{subfigure}[b]{0.45\textwidth}
           \includegraphics[width=1\textwidth]{./figures/coefTM_K8f1.pdf}\vspace{-2mm}
         \caption{$K=8$, $\sigma=1$}
    \end{subfigure}
         \begin{subfigure}[b]{0.45\textwidth}
           \includegraphics[width=1\textwidth]{./figures/coefTM_K8fp2.pdf}\vspace{-2mm}
         \caption{$K=8$, $\sigma=0.2$}
    \end{subfigure}
         \begin{subfigure}[b]{0.45\textwidth}
           \includegraphics[width=1\textwidth]{./figures/coefTM_K2f1.pdf}\vspace{-2mm}
         \caption{$K=2$,  $\sigma=1$}
    \end{subfigure}
         \begin{subfigure}[b]{0.45\textwidth}
           \includegraphics[width=1\textwidth]{./figures/coefTM_K2fp2.pdf}\vspace{-2mm}
         \caption{$K=2$,  $\sigma=0.2$}
    \end{subfigure}
 \fi   
       \caption{Estimated coefficients $(c^{v}_{k,1}, c^{R}_{k,1}, c^{w}_{k,j})$ in NAR models with $p=1$ and $\delta = 5dt$  in four settings of $(K,\sigma)$. The estimators tend to converge fast as the trajectory length $T$ and number $M$ increase: note that the coefficients $c^{w}_{k,1}$ are at the scale of $10^{-4}$ or $10^{-3}$. 
       } \label{fig:coefTM}
\end{figure}

{\bf Consistency of estimators.} The estimator of the NAR models tends to converge as data size increases. Figure \ref{fig:coefTM} shows that the estimated coefficients of NAR with $p=1$ from data consisting of $M$ trajectories, each with length $T$, where $M\in \{2, 8,  32, 128, 512 \}$ and $T\in \{40,80,160,320,640,1280\}$. As $T\times M$ increases, all the estimators tend to converge (note that the coefficients $c^{w}_{k,1}$ are at the scale of $10^{-4}$ or $10^{-3}$). In particular, they converge faster when $\sigma=1$ than when $\sigma=0.2$: the estimators in (a)-(c) oscillate little after $T\times M>10^3$, indicating that different trajectories lead to similar estimators, while the estimators (take $c^R_{K,1}$ for example) in (b)--(d)  oscillate until $T\times M>10^5$. This agrees with the fact that a larger stochastic force makes the system mix faster, so each trajectory provides more effective samples driving the estimator to converge faster.

Numerical tests also show that an NAR model can be numerically unstable, while its coefficient estimator was consistent (i.e., tending to converge as above). Thus, consistency is not sufficient for the selection of an NAR model. 

In our tests, sparse regression algorithms such as LASSO (see e.g., \cite{tibshirani1996_RegressionShrinkage}) or sequential thresholding (see e.g., \cite{she2009_ThresholdingbasedIterative,quade2018_SparseIdentification}) have difficulty in proper thresholding, because the coefficient $c^w$ of the high order terms are small and can vary in scales in different settings, but these high order terms are  important for the NAR model. 

\medskip
Since the NAR models with $p=1$ perform well in all the four settings, and since they are the simplest, we use them in the next sections to explore the maximal time reduction. 

\begin{figure}[H]
    \centering
    \hspace{1mm}
     \ifarXiv
     \includegraphics[width=0.92\textwidth]{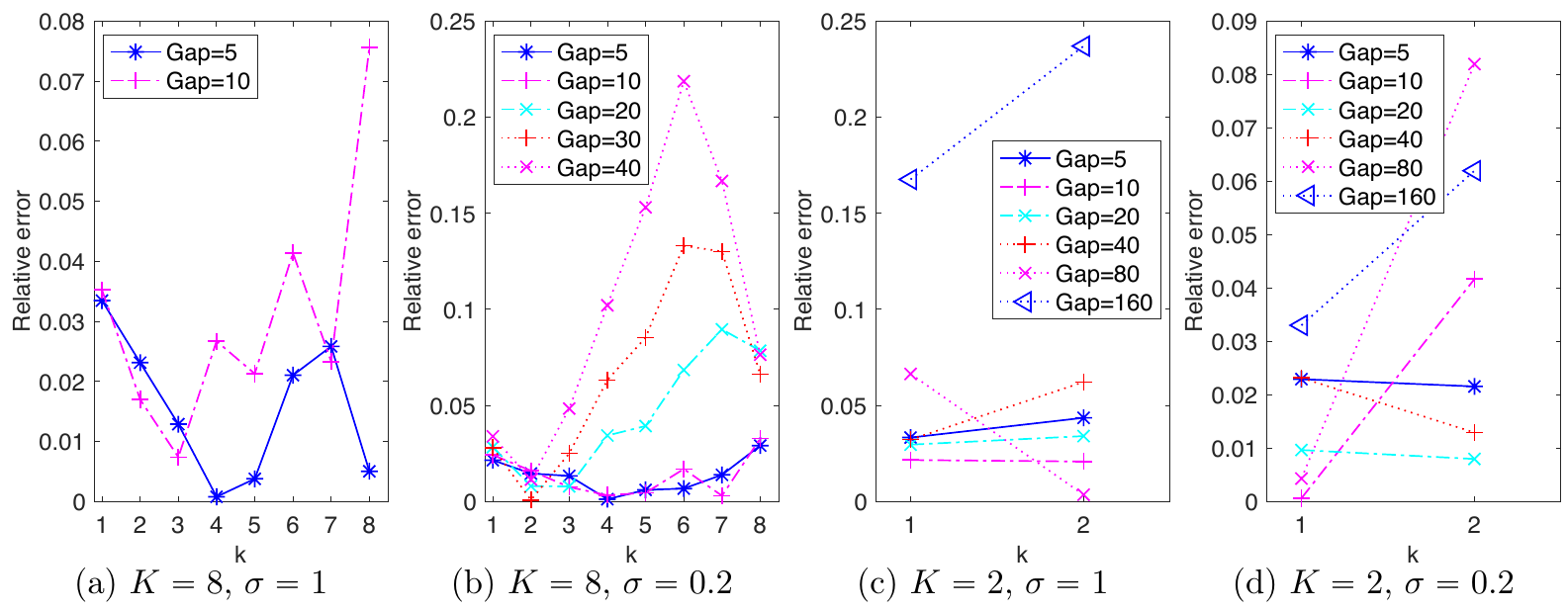} 
\fi
 \ifjournal
     \begin{subfigure}[b]{0.225\textwidth}
           \includegraphics[width=1\textwidth]{./figures/vSpectGaps_Lag1K8f1.pdf}\vspace{-2mm}
         \caption{$K=8$, $\sigma=1$}
    \end{subfigure}
         \begin{subfigure}[b]{0.225\textwidth}
           \includegraphics[width=1\textwidth]{./figures/vSpectGaps_Lag1K8fp2.pdf}\vspace{-2mm}
         \caption{$K=8$, $\sigma=0.2$}
    \end{subfigure}
         \begin{subfigure}[b]{0.225\textwidth}
           \includegraphics[width=1\textwidth]{./figures/vSpectGaps_Lag1K2f1.pdf}\vspace{-2mm}
         \caption{$K=2$,  $\sigma=1$}
    \end{subfigure}
         \begin{subfigure}[b]{0.225\textwidth}
           \includegraphics[width=1\textwidth]{./figures/vSpectGaps_Lag1K2fp2.pdf}\vspace{-2mm}
         \caption{$K=2$,  $\sigma=0.2$}
    \end{subfigure}
\fi    
       \caption{Relative error in energy spectrum reproduced by the NAR models with time steps $\delta= dt\times \rm Gap$ for $\rm Gap\in \{5,10,20,30,40,50\}$ in four settings of $(K,\sigma)$. All NAR models are with time lag $p=1$. The missing $\rm Gap$s in (\textbf{a})--(\textbf{b}) lead to numerically unstable NAR models. Thus, the maximal $\delta$s that an NAR model can reach are $\delta \in [0.01,0.02)$ and $\delta \in [0.04,0.05)$ for (\textbf{a}) and (\textbf{b})  respectively, and $\delta\geq 0.16$ for (\textbf{c})--(\textbf{d}).}  \label{fig:spectrumGaps}
\end{figure}

\subsection{Reduction of the Deterministic Response}\label{sec:K=8}
We explore in this and the next section the maximal time step $\delta$ that the NAR models can reach. We consider only the simplest models with time lag $p=1$. 

We consider first the models with $K=8$ Fourier modes. Since the stochastic force acts directly only on the first $K_0=4$ Fourier modes, the unresolved variable $w$ in \eqref{eq:w_intg} is a deterministic functional of the path of the $K$ modes, so is the truncation error $PB(v+w)- PB(v)$ in \eqref{eq:QBvLag}. Thus, the NAR model mainly reduces the deterministic response of the resolved variables to the unresolved variables. In particular, the term $\Phi^n$ in the NAR model \eqref{eq:NAR_v} optimally approximates this deterministic response on the function space linearly spanned by the terms in \eqref{eq:NAR_linearQ}. 

We consider time steps $\delta = dt\times \gap$ with $\gap\in \{5,10,20,30,40,50\}$. For each $\delta$, we first estimate the coefficients  $(c^{v}_{k,1}, \, c^{R}_{k,1},\, c^{w}_{k,1})$ of the NAR model from the data with the same time step. We then validate the estimated NAR model by its statistics.

Numerical tests show that the NAR models with $\gap \geq 20$ are numerically unstable for the setting $(K=8,\sigma=1)$, and the number is $\gap=50$ for the setting $(K=8,\sigma=0.2)$. Figure \ref{fig:spectrumGaps}a--b shows the relative error in energy spectrum reproduced by NAR models with those stable time steps. The relative errors increase as the $\gap$ increases. Note that the relative errors for modes $k=1,2$ change little, but those with $k\in\{3,4,5,6\}$ increase significantly. In particular, note that in (b), the relative errors at $k=8$ are about 8\% for $\gap\in \{20,30,40\}$, but the relative errors at $k\in\{3,4,5,6\}$ increase sharply to form a peak at $k=6$ when $\gap=40$. We will discuss connections with CFL numbers in Section \ref{sec:CFL}.

\begin{figure}[H]
    \centering
    \hspace{1mm}
 \ifarXiv
     \includegraphics[width=0.92\textwidth]{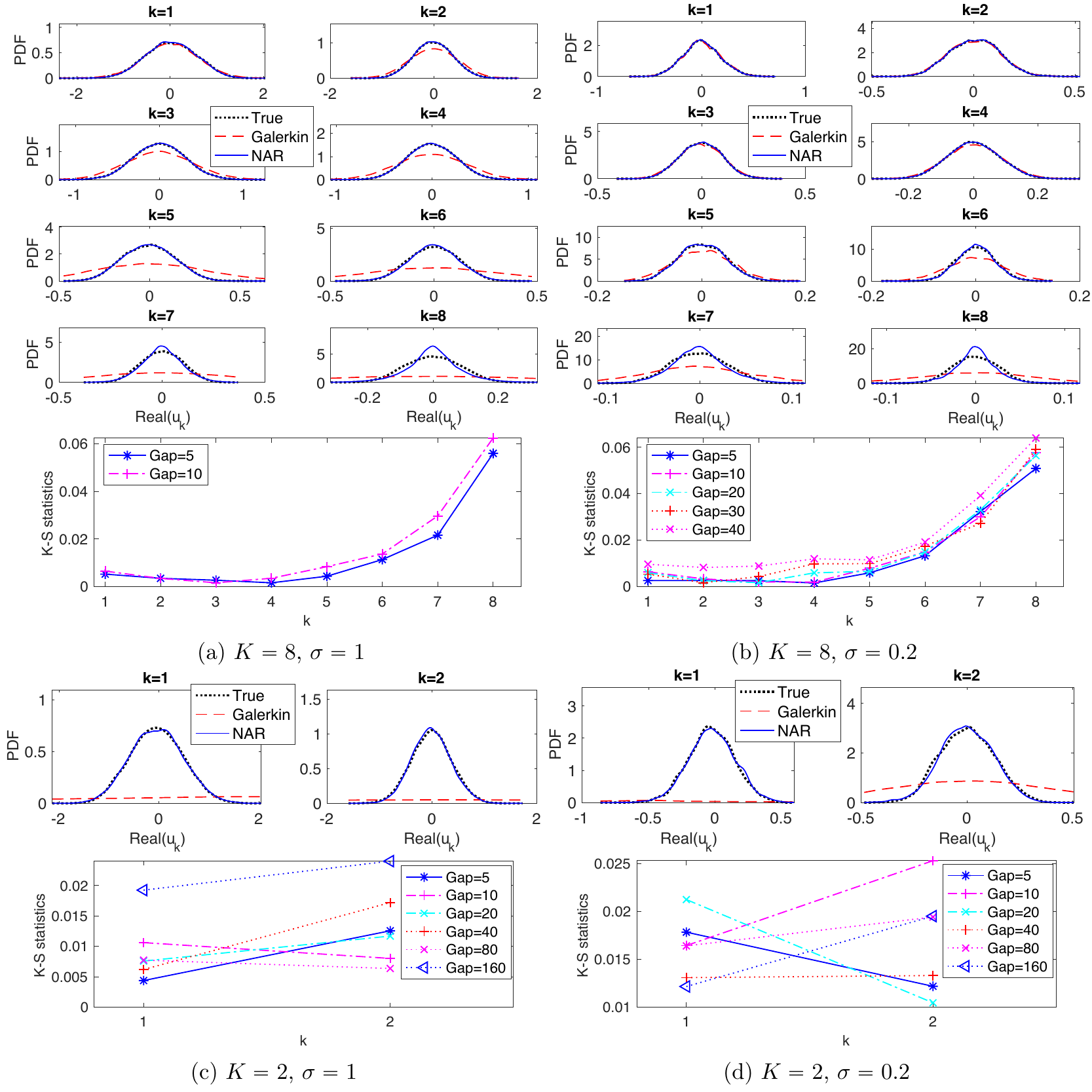} 
\fi
 \ifjournal
       \begin{subfigure}[b]{0.48\textwidth}\centering
     \includegraphics[width=\textwidth]{./figures/pdf_K8f1Lag1.pdf} \\
      \includegraphics[width=0.9\textwidth]{./figures/pdfGaps_K8f1Lag1.pdf}
       \caption{$K=8$, $\sigma=1$}
         \end{subfigure}
              \begin{subfigure}[b]{0.48\textwidth}\centering
     \includegraphics[width=\textwidth]{./figures/pdf_K8fp2Lag1.pdf} 
          \includegraphics[width=0.9\textwidth]{./figures/pdfGaps_K8fp2Lag1.pdf}
      \caption{$K=8$, $\sigma=0.2$}
         \end{subfigure}
                \begin{subfigure}[b]{0.48\textwidth}\centering
     \includegraphics[width=\textwidth]{./figures/pdf_K2f1Lag1.pdf} \\
      \includegraphics[width=0.9\textwidth]{./figures/pdfGaps_K2f1Lag1.pdf}
       \caption{$K=2$, $\sigma=1$}
         \end{subfigure}
              \begin{subfigure}[b]{0.48\textwidth}\centering
     \includegraphics[width=\textwidth]{./figures/pdf_K2fp2Lag1.pdf} 
          \includegraphics[width=0.9\textwidth]{./figures/pdfGaps_K2fp2Lag1.pdf}
      \caption{$K=2$, $\sigma=0.2$}
         \end{subfigure}  
\fi
      \caption{Marginal PDFs and K-S statistics (Kolmogorov--Smirnov statistics, which is the maximum difference between the cumulative distribution functions). In each of (\textbf{a})--(\textbf{d}), the top panels are plots of the empirical marginal PDFs of the real parts of the Fourier modes, from data (True), the K-mode Galerkin system (Galerklin) and the NAR models with $p=1$ and $\delta =\mathrm{Gap} dt$ with $\mathrm{Gap}=5$; the bottom panels are the K-S statistics of NAR models with different time steps $\delta = dt\times {\rm Gap}$, up to the largest $\mathrm{Gap}$ such that the NAR model is numerically stable. } \label{fig:PDFK8}
\end{figure}

These NAR models reproduce the PDFs and ACFs relatively accurately. Figure \ref{fig:PDFK8} shows the marginal PDFs of the real parts of the modes. The top row shows the marginal PDFs for the NAR models with $\gap=5$, in comparison with those of the full model and the Galerkin truncated system (solved with time step $dt$). For the modes with wave numbers $k\in \{1,2,3,4\}$, the NAR model captures the shape and spread of the PDFs almost perfectly, improving those of the Galerkin truncated system. For the modes with $k\in \{5,6,7,8\}$, the NAR model still performs well, significantly improving those of the truncated Galerkin system. The discrepancy between the PDFs becomes larger as the wavenumber increases, because these modes are affected more by the unresolved modes. The bottom row shows that the Kolmogorov--Smirnov statistics (the maximal difference between the cumulative distribution functions) increase slightly as the $\gap$ increases.  

\begin{figure}[H]
    \centering
   \ifarXiv
     \includegraphics[width=0.92\textwidth]{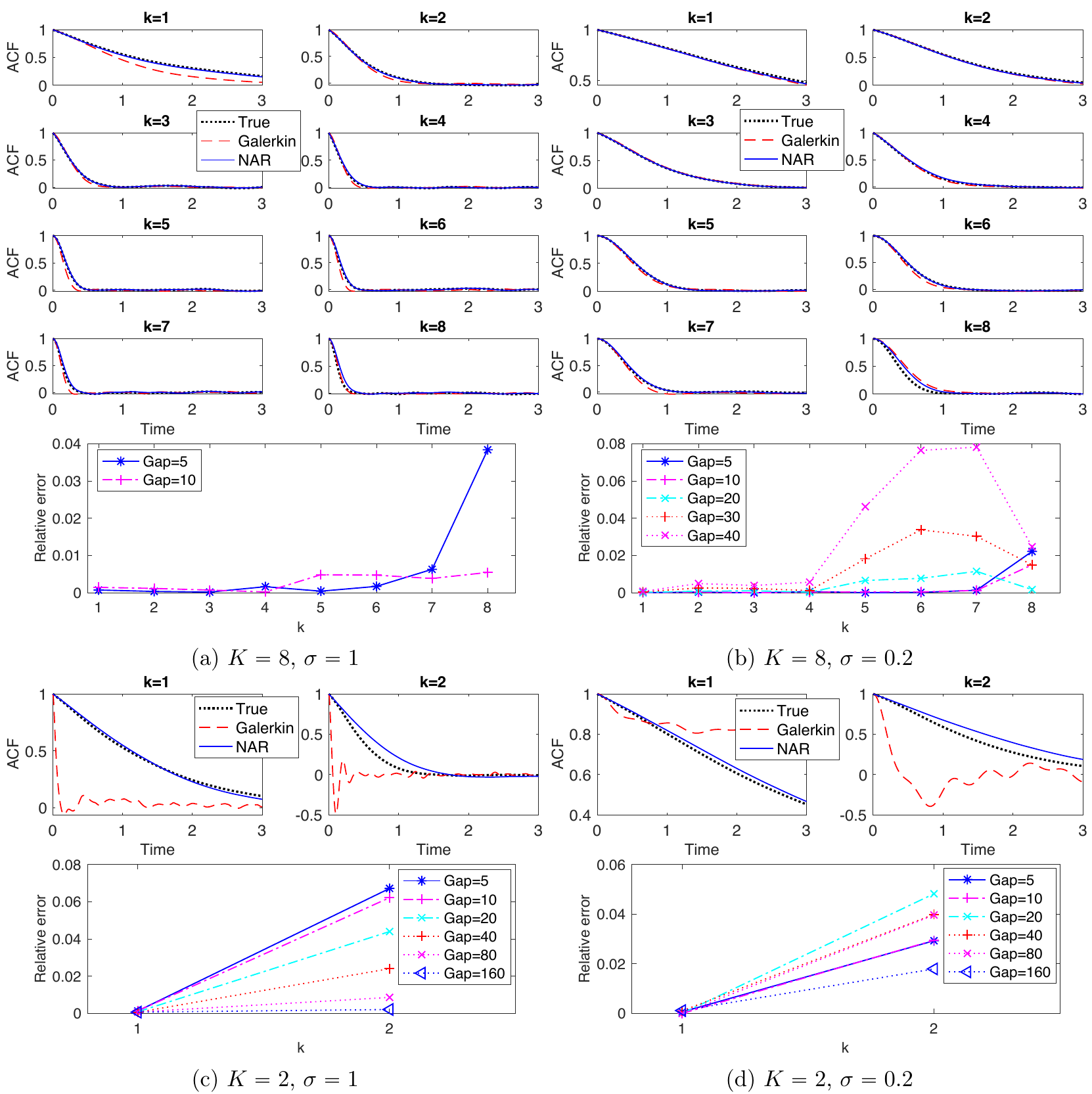} \\
   \fi
   \ifjournal
        \begin{subfigure}[b]{0.48\textwidth}\centering
     \includegraphics[width=\textwidth]{./figures/acf_K8f1Lag1.pdf} \\
      \includegraphics[width=0.9\textwidth]{./figures/acfGaps_K8f1Lag1.pdf}
       \caption{$K=8$, $\sigma=1$}
         \end{subfigure}
              \begin{subfigure}[b]{0.48\textwidth}\centering
     \includegraphics[width=\textwidth]{./figures/acf_K8fp2Lag1.pdf} 
          \includegraphics[width=0.9\textwidth]{./figures/acfGaps_K8fp2Lag1.pdf}
      \caption{$K=8$, $\sigma=0.2$}
         \end{subfigure}
               \begin{subfigure}[b]{0.48\textwidth}\centering
     \includegraphics[width=\textwidth]{./figures/acf_K2f1Lag1.pdf} \\
      \includegraphics[width=0.9\textwidth]{./figures/acfGaps_K2f1Lag1.pdf}
       \caption{$K=2$, $\sigma=1$}
         \end{subfigure}
              \begin{subfigure}[b]{0.48\textwidth}\centering
     \includegraphics[width=\textwidth]{./figures/acf_K2fp2Lag1.pdf} 
          \includegraphics[width=0.9\textwidth]{./figures/acfGaps_K2fp2Lag1.pdf}
      \caption{$K=2$, $\sigma=0.2$}
         \end{subfigure}
         \fi
    \caption{ACF (auto correlation functions). In each of  (\textbf{a})--(\textbf{d}), the top panels are the ACFs of the real parts of the Fourier modes when $\mathrm{Gap}=5$; the  bottom panels are the relative errors (in $L^2([0,3])$-norm) of the NAR models with different time steps $\delta = dt\times {\rm Gap}$, up to the largest $\mathrm{Gap}$ such that the NAR model is numerically stable.} \label{fig:ACFK8}
\end{figure}

Figure \ref{fig:ACFK8} shows the ACFs. The top row shows that both the NAR model (with $\gap=5$) and the Galerkin system can reproduce the ACFs accurately. 
The bottom row shows that the relative error of the ACF, in $L^2([0,3])$-norm, increases as $\gap$ increases (particularly in the case $\sigma=0.2$). 
Recall that the truncated Galerkin system produces PDFs with support much wider than the truth for the high modes (see Figure \ref{fig:PDFK8}), and that $R^\delta$ becomes less accurate as $\delta$ increases. Thus, the terms $u$ and $R^\delta(u)$ in the NAR model  \eqref{eq:NAR} preserve the temporal correlation, and the high order term helps to dissipate energy and preserve the invariant measure.

In summary, when $K=8$, the maximal time steps are $\delta \in dt\times [10,20) = [0.01,0.02)$ and $\delta = dt\times [40,50) = [0.04,0.05)$ when $\sigma =1$ and $\sigma =0.2$, respectively, for NAR models with $p=1$. All these NAR models can accurately reproduce the energy spectrum, the invariant measure and the temporal autocorrelation.

\subsection{Reduction Involving Unresolved Stochastic Force}\label{sec:K=2}
We consider next NAR models with $K=2$. In this case, the unresolved variable $w$ in \eqref{eq:w_intg} is a functional of both the path of the $K$ modes and the unresolved stochastic force. Thus, in view of \eqref{eq:QBvLag} and \eqref{PBu-PBv}--\eqref{eq:NAR}, the NAR model quantifies the response of the $K$-modes to both the unresolved Fourier modes and the unresolved stochastic force. 

Note first that $K=2$ is too small for the $K$-mode Galerkin system to meaningfully reproduce any of the statistical or dynamical properties; see Figure \ref{fig:spectrum}c--d for the energy spectrum, Figure \ref{fig:PDFK8}c--d for the marginal PDFs and Figure \ref{fig:ACFK8}c--d  for the ACFs. On the contrary, the NAR models with $\delta = 5dt$, whose term $R^\delta$ comes from the $K$-mode Galerkin, reproduce these statistics accurately. Remarkably, the NAR models remain accurate even when the time step is as large as $\delta = 80dt$, with the K-S statistics being less than $0.025$ in Figure \ref{fig:PDFK8}c--d,  and with the relative error in ACFs less than 6\% in Figure \ref{fig:ACFK8}c--d.

To explore the maximal time step that NAR models can reach, we consider time steps $\delta = dt\times \gap$ with $\gap\in \{5,10,20,40,80,160\}$. Numerical tests show that the NAR models are numerically stable for all of them in both settings of $\sigma=1$ and $\sigma=0.2$. Figure \ref{fig:spectrumGaps}c--d shows the relative error in energy spectrum reproduced by NAR models with these time steps. The relative error first decreases and then increases as $\gap$ increases, reaching the lowest when $\gap=10$ and $\gap=20$ for the settings $\sigma=1$ and $\sigma=0.2$, respectively. In particular, all of these relative errors remain less than 9\%, except when $\gap=160$ in the setting $\sigma=1$.

In summary, when $K=2$, NAR models can tolerate large time steps. The maximal time steps are at least $\delta =dt\times 80 = 0.08$ and $\delta dt\times 160 = 0.16$ when $\sigma =1$ and $\sigma =0.2$, respectively, for the NAR models to reproduce the energy spectrum with relative error less than 9\%.

\subsection{Discussion on Space-Time Reduction}\label{sec:CFL}
Since model reduction aims for space-time reduction, it is natural to consider the maximal reduction in space-time; in other words, the minimum ``spatial'' dimension $K$ and maximum time step $\delta =dt\times \gap$. We have the following observations from the previous sections:  

\begin{figure}[t]
    \centering
     \includegraphics[width=0.5\textwidth]{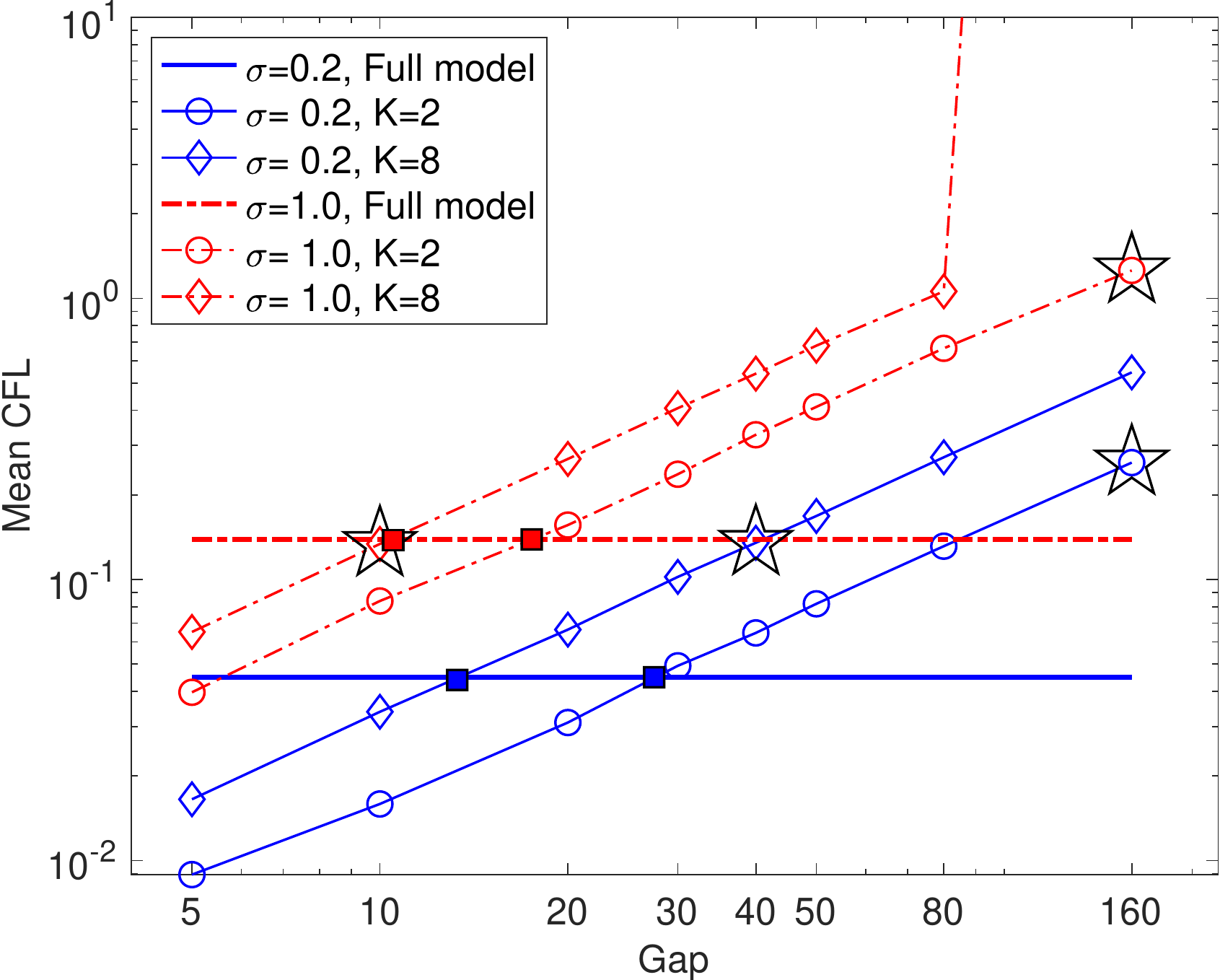} 
      \caption{The mean CFL numbers of the full models and the $K$-mode Galerkin systems. The mean CFL number is computed along a trajectory with $10^5$ steps. The time step is $dt=0.001$ for the full model, and is $\delta=dt\times\gap$ for the $K$-mode Galerkin system. When $(\sigma=1,K=8)$,  the $K$-mode Galerkin system blows up after $\gap>80$, so its CFL number is missing afterwards. The stars are the maximal $\gap$, such that the NAR model is stable. The red and blue squares are where the full model's mean CFL numbers agree with those of the $K$-mode Galerkin systems. The stars ($\medstar$) are the largest time Gap that our NAR model is numerically stable. The relative errors in energy spectrum in Figure \ref{fig:spectrumGaps}c--d are the smallest when the $\gap$'s are the closest to these squares.  } \label{fig:CFL}
\end{figure}

\begin{enumerate}
\item Space dimension reduction, memory length of the reduced model and the stochastic force are closely related. As suggested by the discrete Mori--Zwanzig formalism for random dynamics (see e.g., \cite{LinLu20}), space dimension reduction would lead to non-Markovian closure models. Figure \ref{fig:spectrumLag} suggests that a proper medium length of the memory leads to best NAR model. It also suggests that the scale of the white in time stochastic force can affect the memory length, and a larger scale of stochastic force leads to shorter memory. We leave it as future work to investigate the relations between memory length (colored or white in time), stochastic force, and energy dissipation.  

\item Maximal time step depends on the space dimension and the scale of the stochastic force, mainly limited by the stability of the nonlinear reduced model. Figure \ref{fig:spectrumGaps} shows that the maximum time step when $K=2$ is at least $\delta=dt\times \gap$ with $\gap=160$, much larger than those of the case of $K=8$. It also shows that as the scale of stochastic force increases from $\sigma=0.2$ to $\sigma=1$, the NAR models' maximal time step decreases (because the NAR models either become unstable or have larger errors in energy spectrum). It is noteworthy to mention that these maximal time steps of NAR models are smaller than those that the $K$-mode Galerkin system can tolerate. Figure \ref{fig:CFL} shows that the $K$-mode Galerkin system can be stable for time steps much larger than those of the NAR models: the maximal time step for the K-mode Galerkin system is when the mean CFL number (which increases linearly) reaches 1, but the maximal time step for the NAR models to be stable is smaller. For example, in the setting  $(K=8,\sigma=0.2)$, the maximal time gap for the Galerkin system is $\gap=80$ (the end of the red diamond line), but the maximal time gap for the NAR model is about $\gap=10$.  The increased numerical instability of the NAR model is likely due to the nonlinear terms $\Phi^n$, which are important for the NAR model to preserve energy dissipation and the energy spectrum (see Figure \ref{fig:spectrum} and the coefficients in Figure \ref{fig:coefTM}).  
\end{enumerate}

Beyond maximal reduction, an intriguing question arises: when does the reduced model perform the best (i.e., the least relative error in energy spectrum)? We call it optimality of space-time reduction. It is more interesting and relevant to model reduction than maximal reduction in space-time, because one may achieve a large time step or a small space dimension at the price of a large error in the NAR model, as we have seen in Figure \ref{fig:spectrumGaps}. We note that the relative errors in energy spectrum in Figure \ref{fig:spectrumGaps}c--d are the smallest when the $\gap$s are the closest to the squares in Figure \ref{fig:CFL}, where the full model's mean CFL numbers agree with those of the $K$-mode Galerkin system. We conjecture that optimal space-time reduction can be achieved by an NAR model when the $K$-mode Galerkin system preserves the CFL number of the full model.

\section{Conclusions}
\label{sec5}
We consider a data-driven model reduction for stochastic Burgers equations, casting it as a statistical learning problem on approximating the flow map of low-wavenumber Fourier modes. We derive a class of efficient parametric reduced closure models, based on representing the high modes as functionals of the resolved variables' trajectory. The reduced models are nonlinear autoregression (NAR) time series models, with coefficients estimated from data by least squares. In various settings, the NAR models can accurately reproduce the statistics such as the energy spectrum, the invariant densities, and the autocorrelations. 

Using the simplest NAR model, we investigate the maximal space-time reduction in four settings: reduction of deterministic responses ($K>K_0$) vs.~reduction involving unresolved stochastic force ($K<K_0$), and small vs.~large scales of stochastic force (with $\sigma=0.2$ and $\sigma=1$), where $K_0$ is the number of Fourier modes of the white-in-time stochastic force, and $\sigma$ is the scale of the force. 
 Reduction in space dimension is unlimited, and NAR models with $K=2$ Fourier modes can reproduce the energy spectrum with relative errors less than 5\%. The time reduction is another story. 
 Maximal time reduction depends on both the dimension reduction and the stochastic force's scale, as they affect the stability of the NAR model. The NAR model's stability limits the maximal time step to be smaller than those of the K-mode Galerkin system. Numerical tests indicate that the NAR models achieve the minimal relative error at the time step where the K-mode Galerkin system's mean CFL number agrees with the full model's. This is a potential criterion for optimal space-time reduction. 


The simplicity of our NAR model structure opens various fronts for a further understanding of data-driven model reduction. Future directions include:
(1) studying the connection between optimal space-time reduction, the CFL number, and quantification of the accuracy of reduced models; 
(2) investigating the relation between memory length, dimension reduction, the stochastic force, and the energy dissipation of the system; 
(3) developing post-processing techniques \red{to efficiently recover information of the high Fourier modes}, so as to predict the shocks using the reduced models.


\vspace{6pt}

\paragraph{Acknowledgement} FL is grateful for supports from NSF-1913243 and NSF-1821211. FL would like to thank the anonymous reviewers for valuable feedback that helped significantly improve the manuscript. FL is grateful for Alexandre Chorin for introducing this problem. This study is part of our joint project on renormalization group methods. FL would like to thank Kevin Lin, Panos Stinis, John Harlim, Xiantao Li, Mauro Maggioni, Felix Ye, and Xingjie Li for helpful discussions. 

 \bibliographystyle{plain}

\end{document}